\input amstex
\documentstyle{amsppt}
\document

\magnification 1100

\def\ten{\frak{t}}

\def\sen{\frak{s}}

\def\a{{\alpha}}

\def\ab{{\bold a}}
\def\bb{{\bold b}}
\def\cb{{\bold c}}

\def\eb{{\bold e}}
\def\fb{{\bold f}}
\def\hb{{\bold h}}
\def\ib{{\bold i}}
\def\jb{{\bold j}}
\def\kb{{\bold k}}

\def\lb{{\bold l}}

\def\xb{{\bold x}}

\def\Ab{{\bold A}}
\def\Bb{{\bold B}}

\def\Fb{{\bold F}}
\def\Ib{{\bold I}}
\def\Jb{{\bold J}}

\def\Hb{{\bold H}}

\def\Rb{{\bold R}}

\def\Vb{{\bold V}}

\def\i{{\roman i}}
\def\o{{\roman o}}
\def\d{{\roman d}}
\def\e{{\roman e}}
\def\f{{\roman f}}
\def\h{{\roman h}}

\def\m{{\roman m}}
\def\n{{\roman n}}
\def\n{{\roman n}}
\def\r{{\roman r}}
\def\s{{\roman s}}
\def\t{{\roman t}}
\def\w{{\roman w}}
\def\x{{\roman x}}

\def\E{{\roman E}}

\def\M{{\roman M}}

\def\H{{\roman H}}

\def\A{{\roman A}}

\def\H{{\roman H}}

\def\I{{\roman I}}

\def\K{{\roman K}}

\def\ad{{\roman a\roman d}}

\def\Hom{{\H\o\m\,}}
\def\End{{\E\n\d\,}}
\def\and{{\roman a\roman n\roman d\,}}
\def\Ker{{\K\e\r\,}}
\def\Im{{\I\m\,}}
\def\max{{\m\roman a\x\,}}

\def\where{{\w\h\e\r\e}}
\def\with{{\w\i\t\h}}
\def\or{{\o\r}}

\def\CC{{\Bbb C}}
\def\FF{{\Bbb F}}
\def\KK{{\Bbb K}}

\def\NN{{\Bbb N}}
\def\PP{{\Bbb P}}

\def\ZZ{{\Bbb Z}}

\def\Bc{{\Cal B}}

\def\Oc{{\Cal O}}
\def\Pc{{\Cal P }}

\def\ds{\displaystyle}                
\def\ts{\textstyle}                
\def\ss{\scriptstyle}                
\def\qed{\hfill $\sqcap \hskip-6.5pt \sqcup$}        
\overfullrule=0pt                                    

\def\sln{{{\frak{sl}}_n}}
\def\gln{{{\frak{gl}}_n}}
\def\slh{{\widehat{\frak{sl}}_n}}

\def\sltor{{\ddot{\frak{sl}}_n}}
\def\sltord{{\ddot{\frak{sl}}_{n,d}}}
\def\sltorl{{\ddot{\frak{sl}}_{n,\partial}}}

\def\Uaf{{\dot{\bold U}}}
\def\Uafh{{\dot{\bold U}_{h}}}
\def\Uafv{{\dot{\bold U}_{v}}}
\def\Haf{\dot{\bold H}_m}

\def\Utor{\ddot{\bold U}}
\def\Htor{\ddot{\bold H}_m}

\def\ue{{\underline e}}

\def\oe{{\overline e}}
\def\um{{\underline m}}
\def\ui{{\underline i}}
\def\ul{{\underline l}}
\def\om{{\overline m}}
\def\uk{{\underline k}}
\def\ok{{\overline k}}
\def\uj{{\underline j}}
\def\oi{{\overline i}}
\def\oj{{\overline j}}
\def\uib{{\underline\ib}}
\def\ujb{{\underline\jb}}
\def\ueb{{\underline\eb}}
\def\ojb{{\overline\jb}}

\def\ulb{{\underline\lb}}


\centerline{{\bf DOUBLE-LOOP ALGEBRAS AND}}
\centerline{{\bf THE FOCK SPACE}}

\vskip 12mm

\centerline {\bf M.Varagnolo and E.Vasserot}

\vskip3cm 

\noindent{\bf Introduction.} 
The main motivation of this article comes from physic :
it is related to the Yangian symmetry in conformal field theory
and the spinons basis. In a few words, it has been recently noticed that
level one representations of the affine Lie algebra $\slh$
admit an action of a quantum group, the Yangian of type $A_n^{(1)}$.
A quantized version of this statement says that
the Fermionic Fock space admits
two different actions of the quantized enveloping algebra of $\slh$.
The first one is a $q$-deformation of the
well-known level-one representation of the affine Lie algebra
(see [H], [KMS]).
When $q$ is one this representation may be viewed as a
particular case of the Borel-Weil theorem for loop groups (see [PS]). 
The second one is a level-zero action
arising from solvable lattices models (more precisely the
Calogero-Sutherland and the Haldane-Shastry models,
see [JKKMP], [TU], and the references therein).
Quite remarkably these two constructions can be glued together to 
get a representation of a new object (introduced in [GKV] and [VV]) : 
a toroidal quantum group, i.e.
a two parameters deformation of the enveloping algebra of the universal
extension of the Lie algebra $\sln[x^{\pm 1},y^{\pm 1}]$.
The aim of this note is three-fold. 
First we define a representation of the quantized toroidal algebra, $\Utor$, 
on the Fock space 
generalizing the two actions of the affine quantum group previously known.
For that purpose we first construct an action on the space
$\bigwedge^m(\CC^n[z^{\pm 1}])$ for any positive integer $m$ by means
of the Schur-type duality between $\Utor$ and Cherednik's double affine Hecke 
algebra established in [VV], then we explain how to perform the limit 
$m\to\infty$.
The second purpose of this article is to explain to which extend this 
representation can be viewed 
in geometrical terms, by means of correspondences on infinite flags
manifolds.
A complete geometric picture would require equivariant $K$-theory
of some infinite dimensional variety. The correct definition of such
$K$-groups will be done in another work (see [GKV] and [GG] for
related works). We will mainly concentrate here on the algebraic
aspects.
At last, an essential point in the Fock space representation that we consider 
is that it involves some polynomial difference operators. 
It is due to the fact, proved in section 13,
that the classical toroidal algebra, i.e. the 
specialization to $q=1$ of the toroidal algebra, is isomorphic to 
the enveloping algebra of the universal central extension of a 
current Lie algebra over a quantum torus. 

\vskip1mm

\noindent Y. Saito, K. Takemura and D. Uglov have obtained
similar results in [STU]. The computations in the proof
of the formulas (12.7-8) (formula (6.16) in [STU]),
not written in the first version of our paper,
are different but rely on the same results from [TU] and [VV].

\vskip3mm 

\noindent{\bf 1.} Fix $q\in\CC^\times$.
The toroidal Hecke algebra of type ${\frak{gl}}_m$, $\Htor$, is the 
unital
associative algebra over $\CC[\xb^{^{\pm 1}}]$
with the generators 
$T_i^{^{\pm 1}},$ $X_j^{^{\pm 1}},$ $Y_j^{^{\pm 1}},$
$i=1,2,...,m-1,$ $j=1,2,...,m$
and the relations 
$$T_i\,T_i^{-1}=T_i^{-1}\,T_i=1,\qquad (T_i+q^{-1})\,(T_i-q)=0,$$
$$T_i\,T_{i+1}\,T_i=T_{i+1}\,T_i\,T_{i+1},$$
$$T_i\,T_j=T_j\,T_i\quad {\i\f}\ |i-j|>1,$$
$$X_0\,Y_1=\xb\,Y_1\,X_0,
\qquad X_i\,X_j=X_i\,X_j,\qquad Y_i\,Y_j=Y_j\,Y_i,$$
$$X_j\,T_i=T_i\,X_j,\qquad Y_j\,T_i=T_i\,Y_j,
\qquad {\i\f}\ j\not= i,i+1$$
$$T_i\,X_i\,T_i=X_{i+1},\qquad 
T_i^{-1}\,Y_i\,T_i^{-1}=Y_{i+1},$$
$$Y_2\,X_1^{-1}\,Y_2^{-1}\,X_1=
T_1^{-2},$$
where $X_0=X_1\,X_2\cdots X_m$.
The algebra $\Htor$ has been introduced previously by Cherednik
to prove the Macdonald conjectures (see [C1]).
Let us mention that we have fixed the central element $\xb$ 
in a different way than in [C1].
Put $Q=X_1\,T_1\,\cdots T_{m-1}\in\Htor$. Then, 
$T_i^{\pm 1},$ $Y_j^{\pm 1},$ $Q^{\pm 1}$,
$i=1,2,...,m-1,$ $j=1,2,...,m$, is a system of generators of $\Htor$.
Besides, for any $i=1,2,...,m-1$ a direct computation gives 
$Q\,Y_i\,Q^{-1}=Y_{i+1}$, and
$Q\,Y_m\,Q^{-1}=\xb\,Y_1$.
Indeed we have 

\vskip.3cm

\noindent{\bf Proposition 1.} {\it The toroidal Hecke algebra $\Htor$ admits a
presentation in terms of generators $T_i^{^{\pm 1}},$ $Y_j^{^{\pm 1}},$ 
$Q^{^{\pm 1}},$ $i=1,2,...,m-1,$ $j=1,2,...,m,$
with relations } 
$$T_i\,T_i^{-1}=T_i^{-1}\,T_i=1,\qquad (T_i+q^{-1})\, (T_i-q)=0,$$
$$T_i\,T_{i+1}\,T_i=T_{i+1}\,T_i\,T_{i+1},$$
$$T_i\,T_j=T_j\,T_i\quad {if}\ |i-j|>1,$$
$$Y_i\,Y_j=Y_j\,Y_i,\qquad
T_i^{-1}\,Y_i\,T_i^{-1}=Y_{i+1},$$
$$Y_j\,T_i=T_i\,Y_j,
\qquad {if}\ j\not= i,i+1$$
$$Q\,T_{i-1}\,Q^{-1}=T_i\quad (1<i<m),\qquad
Q^2\,T_{m-1}\,Q^{-2}=T_1,$$
$$Q\,Y_i\,Q^{-1}=Y_{i+1}\quad (1\leq
i\leq m-1),\qquad  Q\,Y_m\,Q^{-1}=\xb\,Y_1.$$

\vskip 3mm

\noindent{\bf Remarks.} 
{\bf 1.1.} 
Given a permutation $w\in{\frak{ S}}_m$ and a reduced decomposition 
$w=s_{_{i_1}}s_{_{i_2}}\cdots s_{_{i_k}}$ in terms of the simple transpositions
$s_{_1}$, $s_{_2}$,...,$s_{_{m-1}}$, set as usual 
$T_w=T_{i_1}T_{i_2}\cdots T_{i_k}\in\Htor$.
In particular $T_w$ is independent of the choice of the reduced 
decomposition.
Moreover, given a $m$-tuple of integers $\ab=(a_1,a_2,...,a_m)$, denote
by $X^\ab$ and $Y^\ab$ the corresponding monomials in the $X_i$'s
and $Y_i$'s. It is known that the $X^\ab\,Y^\bb\,T_w$'s form a basis of $\Htor$.

\noindent {\bf 1.2.} 
One consequence of the existence of the 
basis of monomials above is that 
the subalgebra of $\Htor$ generated by the $T_i$'s and the $X_i$'s
is isomorphic to the affine Hecke algebra of type ${\frak{gl}}_m$.
Let denote by $\Haf$ this subalgebra. Similarly, the subalgebra 
$\Hb_m\subset\Htor$ generated by the $T_i$'s alone is isomorphic to the finite
Hecke algebra of type ${\frak{gl}}_m$.

\vskip3mm

\noindent The maps
$$\matrix
\omega\,:\qquad T_i,\,Q,\,Y_i,\,\xb,\,q\,\mapsto\,
-T_i^{-1},\,(-q)^{m-1}Q,\,Y_i^{-1},\,\xb^{-1},\,q,\hfill\cr\cr
\gamma\,:\qquad T_i,\,X_i,\,Y_i,\,\xb,\,q\,\mapsto\,
T_i,\,Y_i^{-1},\,X_i^{-1},\,\xb,q,\hfill
\endmatrix$$
extend uniquely to an involution and an anti-involution of $\Htor$.
Let us remark that, if $S_m,A_m\in\Hb_m$ are the 
symmetrizer and the antisymmetrizer of $\Hb_m$, that is to say
$$S_m=\sum_{w\in{\frak{S}}_m}q^{l(w)}T_w\quad{\and}\quad
A_m=\sum_{w\in{\frak{S}}_m}(-q)^{-l(w)}T_w,$$
where $l\,:\,{\frak{S}}_m\to\NN$ is the length, 
then $\omega(S_m)=q^{m(m-1)}A_m$.

\vskip3mm

\noindent{\bf 2.} Fix $p\in\CC^\times$ and set
$\Rb_m=\CC[z_1^{\pm 1},z_2^{\pm 1},...,z_m^{\pm 1}]$.
Consider the following operators in $\End(\Rb_m)$ :
$$\matrix
{\displaystyle t_{i,j}=(1+s_{i,j})\,{q^{-1}z_i-qz_j\over z_i-z_j}-q^{-1},}\qquad\hfill
&1\leq i,j\leq m,\hfill\cr\cr
x_i=t_{i-1,i}\,s_{i-1,i}\cdots t_{1,i}\,s_{1,i}\,D_i\,
s_{i,m}\,t_{i,m}^{-1}\cdots s_{i,i+1}\,t_{i,i+1}^{-1},\qquad\hfill
&i=1,2,...,m,\hfill\cr\cr
y_i=z_i^{-1},\qquad\hfill
&i=1,2,...,m,\hfill
\endmatrix$$
where $s_{i,j}$ acts on Laurent polynomials by permuting $z_i$ and
$z_j$, and $D_i$ is the difference operator such that
$(D_if)(z_1,z_2,...,z_m)=f(z_1,...,pz_i,...,z_m).$
The following result was first noticed by Cherednik.
 
\vskip3mm

\noindent{\bf Proposition 2.} {\it The map 
$$T_i\mapsto t_{i,i+1},\quad Y_i\mapsto y_i,\quad
Q\mapsto D_1\,s_{1,m}\,s_{1,m-1}\cdots s_{1,2},\quad
\xb\mapsto p,$$ 
extends to a representation of $\Htor$ in $\Rb_m$. Moreover,
in this representation $X_i$ acts as $x_i$.}
\qed

\vskip3mm

\noindent{\bf Remark 2.}  
Let us recall that $\Haf\subset\Htor$ is the subalgebra 
generated by the $T_i$'s and the $X_i$'s. The trivial module of $\Haf$ 
is the one-dimensional representation such that
$T_i$ acts by $q$ and $X_i$ by $q^{2i-m-1}$.
Then $\Rb_m$ is the $\Htor$-module induced from the trivial representation 
of $\Haf$. In other words, if $\Ib_m\subset\Htor$ is the left ideal 
generated by the $T_i-q$'s and the $X_i-q^{2i-m-1}$'s, then $\Rb_m$ is 
identified with the quotient $\Htor/\Ib_m$ where $\Htor$ acts by left 
translations.

\vskip3mm

\noindent{\bf 3.} Fix $d\in\CC^\times$ and an integer $n\geq 3$.
The toroidal quantum group of type 
${\frak{sl}}_n$, $\Utor$, is the complex unital associative algebra 
generated by $\eb_{i,k},$ $\fb_{i,k},$ $\hb_{i,l},$ $\kb_i^{\pm 1}$, 
where $i=0,1,...,n-1,$ $k\in\ZZ$, $l\in\ZZ^\times$,
and the central elements $\cb^{^{\pm 1}}.$
The relations are expressed in term of the formal series
$$\eb_i(z)=\sum_{k\in\ZZ}{\eb_{i,k}\cdot z^{^{-k}}},\quad
\fb_i(z)=\sum_{k\in\ZZ}{\fb_{i,k}\cdot z^{^{-k}}},$$
and 
$\kb_i^{^\pm}(z)=\kb_i^{\pm 1}\cdot\exp\left(\pm(q-q^{-1})
\sum_{k\geq 1}{\hb_{i,\pm k}\cdot z^{^{\mp k}}}\right),$
as follows
$$\kb_i\cdot\kb_i^{-1}=\cb\cdot\cb^{^{-1}}=1,\qquad
[\kb^{^\pm}_i(z),\kb^{^\pm}_j(w)]=0,$$
$$\theta_{-a_{ij}}(\cb^{^{2}}d^{^{-m_{ij}}}w z^{^{-1}})\cdot\kb^{^+}_i(z)\cdot
\kb^{^-}_j(w)=\theta_{-a_{ij}}(\cb^{^{-2}}d^{^{-m_{ij}}}wz^{^{-1}})\cdot
\kb^{^-}_j(w)\cdot\kb^{^+}_i(z),$$
$$\kb^{^\pm}_i(z)\cdot\eb_j(w) =
\theta_{\mp a_{ij}}(\cb^{^{-1}}d^{^{\mp m_{ij}}}w^{^{\pm 1}}z^{^{\mp 1}})
\cdot\eb_j(w)\cdot\kb^{^\pm}_i(z),
$$
$$\kb^{^\pm}_i(z)\cdot\fb_j(w) =
\theta_{\pm a_{ij}}(\cb d^{^{\mp m_{ij}}}w^{^{\pm 1}}z^{^{\mp 1}})
\cdot\fb_j(w)\cdot\kb^{^\pm}_i(z),
$$
$$(q-q^{-1})[\eb_i(z),\fb_j(w)]=\delta(i=j)\,
\biggl(\epsilon(\cb^{^{-2}}\cdot z/w)\cdot\kb^{^+}_i(\cb\cdot w)-
\epsilon(\cb^{^2}\cdot z/w)\cdot\kb^{^-}_i(\cb\cdot z)\biggr),$$
$$(d^{^{m_{ij}}}z-q^{a_{ij}}w)\cdot\eb_i(z)\cdot\eb_j(w)=
(q^{a_{ij}}d^{^{m_{ij}}}z-w)
\cdot\eb_j(w)\cdot\eb_i(z),$$
$$(q^{a_{ij}}d^{^{m_{ij}}}z-w)\cdot\fb_i(z)\cdot\fb_j(w)=
(d^{^{m_{ij}}}z-q^{a_{ij}}w)
\cdot\fb_j(w)\cdot\fb_i(z),$$
$$\{\eb_i(z_{_1})\cdot\eb_i(z_{_2})\cdot\eb_{j}(w)-
(q+q^{-1})\cdot\eb_i(z_{_1})\cdot\eb_{j}(w)\cdot\eb_i(z_{_{2}})+
\eb_{j}(w)\cdot\eb_i(z_{_{1}})\cdot\eb_i(z_{_2})\}+$$
$$+\{z_{_1}\leftrightarrow z_{_2}\}=0,
\qquad{\i\f}\quad a_{ij}=-1,$$
$$\{\fb_i(z_{_1})\cdot\fb_i(z_{_2})\cdot\fb_{j}(w)-
(q+q^{-1})\cdot\fb_i(z_{_1})\cdot\fb_{j}(w)\cdot\fb_i(z_{_{2}})+
\fb_{j}(w)\cdot\fb_i(z_{_{1}})\cdot\fb_i(z_{_2})\}+$$
$$+\{z_{_1}\leftrightarrow z_{_2}\}=0,
\qquad{\i\f}\quad a_{ij}=-1,$$
$$[\eb_i(z),\eb_j(w)]=[\fb_i(z),\fb_j(w)]=0\qquad
{\i\f}\quad a_{ij}=0,$$

\vskip2mm

\noindent where $\epsilon(z) = \sum_{_{n=-\infty}}^{^\infty}z^n$,
$\theta_m(z)\in\CC[[z]]$ is the expansion of
${q^m\cdot z -1\over z-q^m}$, and
$a_{ij}$, $m_{ij}$, are the entries of the following $n\times n$-matrices 
\hfill\break

$\A=\pmatrix
2&-1&&0&-1\cr
-1&2&\cdots&0&0\cr
&\vdots&\ddots&\vdots&\cr
0&0&\cdots&2&-1\cr
-1&0&&-1&2
\endpmatrix,
\qquad
\M=\pmatrix
0&-1&&0&1\cr
1&0&\cdots&0&0\cr
&\vdots&\ddots&\vdots&\cr
0&0&\cdots&0&-1\cr
-1&0&&1&0
\endpmatrix.$
\hfill\break

\vskip3mm

\noindent 
Let $\Uaf$ be the quantized enveloping algebra of
$\slh$, i.e. the algebra generated by $\eb_i,\fb_i,\kb_i^{\pm 1}$
with $i=0,1,...,n-1$ modulo the Kac-Moody type relations

$$\kb_i\cdot\kb_i^{\pm 1}=1,\qquad \kb_i\cdot\kb_j=\kb_j\cdot\kb_i,$$

$$\kb_i\cdot\eb_j=q^{a_{ij}}\,\eb_j\cdot\kb_i,\qquad
\kb_i\cdot\fb_j=q^{-a_{ij}}\,\fb_j\cdot\kb_i,$$

$$[\eb_i,\fb_j]=\delta(i=j){\kb_i-\kb_i^{-1}\over q-q^{-1}},$$

\noindent and, if $i\neq j$,

$$\sum_{k=0}^{1-a_{ij}}(-1)^k\eb_i^{(k)}\eb_j\eb_i^{(1-a_{ij}-k)}=
\sum_{k=0}^{1-a_{ij}}(-1)^k\fb_i^{(k)}\fb_j\fb_i^{(1-a_{ij}-k)}=0,$$
where

$$
\matrix
\eb_i^{(k)}=\eb_i^k/[k]!,\qquad &\fb_i^{(k)}=\fb_i^k/[k]!,\cr\cr
[k]={q^k-q^{-k}\over q-q^{-1}},\qquad &[k]!=[k][k-1]\cdots[1].
\endmatrix
$$
\vskip2mm
\noindent Let us recall that $\Uaf$ admits another presentation, the
Drinfeld new presentation (see [D]), similar to the presentation
of $\Utor$ above. The isomorphism between the two presentations
of $\Uaf$ is announced in [D] and proved in [B].

\vskip2mm

\noindent As indicated in [GKV] the
algebra $\Utor$ contains two remarkable subalgebras, 
$\Uafh$ and $\Uafv$, both
isomorphic to a quotient of $\Uaf$. 
The first one, the horizontal subalgebra, is generated by 
$\eb_{i,0}$, $\fb_{i,0}$, $\kb_i^{\pm 1}$, with
$i=0,1,...,n-1$. These elements satisfy the above relations. 
The second one, the vertical subalgebra, is generated by 
$d^{ik}\eb_{i,k},$ $d^{ik}\fb_{i,k},$ $d^{il}\hb_{i,l},$ $\kb_i^{\pm 1}$, 
where $i=1,2,...,n-1,$ $k\in\ZZ$ and $l\in\ZZ^\times$.
These elements satisfy the relations of the new presentation
of $\Uaf$. Fix $\eb_i=\eb_{i,0}$ and $\fb_i=\fb_{i,0}$, for any
$i=0,1,...,n-1$. It is convenient
to fix an additional triple of elements $\eb_n$, $\fb_n$ 
and $\kb_n^{\pm 1}$ such that $\Uafv$ is generated by
$\eb_i$, $\fb_i$, $\kb_i^{\pm 1}$, with $i=1,2,...,n$, 
satisfying the previous Kac-Moody type relations.

\vskip 3mm

\noindent{\bf 4.} 
For any complex vector space $V$ and any formal variable $\zeta$
denote by $V[\zeta^{\pm 1}]$ the tensor product $V\otimes\CC[\zeta^{\pm 1}].$
Let $v_1,v_2,...,v_n$ be a basis of $\CC^n$.
Set $v_{i+nk}=v_i\,\zeta^{-k}$ for all $i=1,2,...,n$ and all $k\in\ZZ$.
The vectors $v_i$, $i\in\ZZ$, form a basis of $\CC^n[\zeta^{\pm 1}]$.
Given $k\in\ZZ$, write $k=n\,\uk+\ok$, where $\uk$ is a certain integer
and $\ok\in\{1,2,...,n\}$. 
The space $\CC^n[\zeta^{\pm 1}]$ is endowed with a 
representation of the quantized enveloping algebra
of $\slh$, $\Uaf$, such that the Kac-Moody generators act as 
$$\matrix
\eb_i(v_j)=\delta(\oj=\overline{i+1}) v_{j-1},\hfill\cr\cr
\fb_i(v_j)=\delta(\oj=\oi) v_{j+1},\hfill\cr\cr
\kb_i(v_j)=q^{\delta(\oj=\oi)-\delta(\oj=\overline{i+1})}\,v_j,\hfill
\endmatrix$$
for all $j\in\ZZ$ and $i=0,1,2...,n-1$, where $\delta(P)$ is 
1 if the statement $P$ is true and 0 otherwise.

\vskip3mm

\noindent{\bf 5.} Consider now the tensor product 
$\bigotimes^m\CC^n[\zeta^{\pm 1}]=
(\CC^n)^{\otimes m}[\zeta_1^{\pm 1},...,\zeta_m^{\pm 1}]$.
The monomials in the $v_i$'s are parametrized by $m$-tuple of integers.
Such a $m$-tuple can be viewed as a function $\jb:\,\ZZ\to\ZZ$ such that
$\jb(k+m)=\jb(k)+n$ for all $k$ : the map $\jb$ is simply
identified with the $m$-tuple 
$(j_1,j_2,...,j_m)=(\jb(1),\jb(2),...,\jb(m))$.
Let $\Pc_m$ be the set of all such functions.
Then, $\bigotimes^m\CC^n[\zeta^{\pm 1}]$
is endowed with a $\Uaf$ action generalizing the representation in 
section 4 as follows (see [GRV]) : if $\jb\in\Pc_m$,
$$\matrix
{\ds \eb_i(v_\jb)=q^{-\sharp\jb^{-1}(i)}\sum_{k\in\jb^{-1}(i+1)}
q^{2\sharp\{l\in\jb^{-1}(i)\,|\,l>k\}}v_{\jb_k^-},}\hfill\cr\cr
{\ds \fb_i(v_\jb)=q^{-\sharp\jb^{-1}(i+1)}\sum_{k\in\jb^{-1}(i)}
q^{2\sharp\{l\in\jb^{-1}(i+1)\,|\,l<k\}}v_{\jb_k^+},}\hfill\cr\cr
{\ds \kb_i(v_\jb)=q^{\sharp\jb^{-1}(i)-\sharp\jb^{-1}(i+1)}\,v_\jb,}\hfill
\endmatrix$$
where $i=0,1,...,n-1,$ 
$v_\jb=v_{j_1}\otimes v_{j_2}\otimes\cdots\otimes v_{j_m}$
and $\jb^\pm_k$ is the function associated to the
$m$-tuple $(j_1,j_2,...,j_k\pm 1,...,j_m)$.
This action commutes with the action of $\Haf$ such that :

\vskip2mm

\item{\bf .} $T_k$, $k=1,2,...,m-1$, is represented by 
$\tau_{k,k+1}$, where $\tau_{k,l}$ is the automorphism of
$\bigotimes^m\CC^n[\zeta^{\pm 1}]$ which  
acts on the $k$-th and $l$-th components as, $\forall\,i,j\in\ZZ$,
$$v_{ij}=v_i\otimes v_j\mapsto\cases
q\,v_{ij}&\text{\quad if \quad $i=j,$}\cr
q^{-1}v_{ji}&\text{\quad if \quad $i<j,$}\cr
q\,v_{ji}+(q-q^{-1})v_{ij}&\text{\quad if \quad $i>j,$}\cr
\endcases$$
and which acts trivially on the other components,

\vskip1mm

\item{\bf .} $Q$ is represented by 
$\vartheta\,:\,v_{i_1\,i_2...i_m}\mapsto
v_{i_m-n,i_1,...,i_{m-1}}=v_{i_m\,i_1...i_{m-1}}\,\zeta_1.$

\vskip2mm

\noindent We do not prove here that the operators above
satisfy the relations of $\Uaf$ and $\Haf$ since it will follow immediatly
from the results in section 7 or the results in section 6.

\vskip3mm

\noindent{\bf 6.} Geometrically the formulas in the previous section may be
viewed as follows. Suppose that $q$ is a prime power
and let $\FF$ be the field with $q^2$ elements.
Denote by $\KK=\FF((z))$ the field of Laurent power series.
A lattice in $\KK^m$ is a free $\FF[[z]]$-submodule of $\KK^m$ of rank $m$.
Let $\Bc$ be the set of complete periodic flags, i.e. of
sequences of lattices $L=(L_i)_{i\in\ZZ}$ such that
$$L_i\subset L_{i+1},\qquad\dim_{_\FF}(L_{i+1}/L_i)=1\quad{\and}\quad
L_{i+m}=L_i\cdot z^{-1}.$$
Similarly $\Bc^n$ is the set of $n$-steps periodic flags in $\KK^m$, i.e. of
sequences of lattices $L^n=(L^n_i)_{i\in\ZZ}$ such that
$$L^n_i\subseteq L^n_{i+1},\quad{\and}
\quad L^n_{i+n}=L^n_i\cdot z^{-1}.$$
The group $GL_m(\KK)$ acts on $\Bc$ and $\Bc^n$ in a natural way.
The orbits of the diagonal action of $GL_m(\KK)$ on $\Bc^n\times\Bc$ are 
parametrized by $\Pc_m$.
If $e_1,e_2,...,e_m\in\KK^m$ is the canonical basis and $e_{i+mk}=e_i\,z^{-k}$
for all $k\in\ZZ$, we associate
to $\jb\in\Pc_m$ the orbit, say $\Oc_\jb$, of the pair $(L^n_\jb,L)$ where
$$L^n_{\jb,i}=\prod_{\jb(j)\leq i}\FF\,e_j\quad{\and}\quad
L_i=\prod_{j\leq i}\FF\,e_j.$$
Let $\CC_{GL_m(\KK)}[\Bc\times\Bc]$ and $\CC_{GL_m(\KK)}[\Bc^n\times\Bc^n]$ 
be the convolution algebras of invariant complex functions 
supported on a finite number of orbits. It is well known that 
$\CC_{GL_m(\KK)}[\Bc\times\Bc]$ is isomorphic to $\Haf$ (see [IM]). 
Similarly if $i=0,1,...,n-1$ let 
$m_i,\chi^\pm_i,\chi^0\in\CC_{GL_m(\KK)}[\Bc^n\times\Bc^n]$ be such that 

\vskip2mm

\itemitem{\bf .} $m_i(L',L)=\dim(L_i/L_0),\qquad\forall L,L'\in\Bc^n$,

\vskip2mm

\itemitem{\bf .} $\chi^\pm_i$ is the characteristic function of the set
$$\{(L^\pm,L^\mp)\in\Bc^n\times\Bc^n\,|\,L_j^-\subset L_j^+\quad\and
\quad \dim(L_j^+/L_j^-)=\delta(\oi=\oj),\quad\forall j\in\ZZ\},$$

\vskip2mm

\itemitem{\bf .} $\chi^0$ is the characteristic function of the diagonal 
in $\Bc^n\times\Bc^n$.

\vskip2mm

\noindent Then the map
$$\eb_i\mapsto q^{m_{i-1}-m_i}\chi^+_i,\quad
\fb_i\mapsto q^{m_i-m_{i+1}}\chi^-_i,\quad
\kb_i\mapsto q^{2m_i-m_{i-1}-m_{i+1}}\chi^0,$$
extends to an algebra homomorphism 
$\Uaf\to\CC_{GL_m(\KK)}[\Bc^n\times\Bc^n]$.
This statement is stated without a proof in [GV; Theorem 9.2].
For the convenience of the reader a proof is given in the appendix.
Let us mention however that this computation is nothing but an adaptation
of the non affine case proved in [BLM].
As a consequence, the convolution product induces an action of 
$\Uaf$ and $\Haf$ on $\CC_{GL_m(\KK)}[\Bc^n\times\Bc]$.

\vskip3mm

\noindent{\bf Proposition 6.} {\it The isomorphism of vector spaces
$\bigotimes^m\CC^n[\zeta^{\pm 1}]
\,{\buildrel\sim\over\longrightarrow}\,
\CC_{GL_m(\KK)}[\Bc^n\times\Bc]$
mapping $v_\jb$ to the characteristic function of the orbit $\Oc_\jb$ 
is an isomorphism of $\Uaf\times\Haf$-modules between the representation
of section 5 and the representation by convolution.
}

\vskip3mm

\noindent{\it Proof.}
Put $G=GL_m(\KK)$ and let $B\subset G$ be the Iwahori subgroup, i.e. 
the subgroup
of matrices mapping each $e_i$ to a linear combination of the type
$$\sum_{j\leq i}a_{ij}\,e_j\quad{\with}\quad a_{ij}\in\FF,\quad
a_{ii}\neq 0.$$
Let us first compute the $\Haf$-action.
In order to simplify the notations we fix $m=2$.
For any $a,b\in\ZZ$ let $L(a,b)\subset\KK^2$ be the lattice
with basis $(e_{1+2a}, e_{2+2b})$, i.e.
$$L(a,b)=\left(\prod_{k\leq a}\FF\, e_{1+2k}\right)\oplus
\left(\prod_{k\leq b}\FF\, e_{2+2k}\right).$$
Let us recall that the element $L\in\Bc$ is the sequence of lattices such that
$L_0=L(-1,-1)$, $L_1=L(0,-1)$ and $L_{i+2}=L_i\cdot z^{-1}$
for any $i\in\ZZ$.
In particular,
$$\Oc_\jb\cap(\Bc^n\times\{L\})=(B\cdot L_\jb^n)\times\{L\}.$$
By definition, the isomorphism 
$\Haf{\buildrel\sim\over\rightarrow}\CC_G[\Bc\times\Bc]$
maps $T_1$ to $q^{-1}$ times the characteristic function of the $G$-orbit 
$$G\cdot (L',L)\subset\Bc\times\Bc,$$
where $L'_0=L_0$ and $L'_1=L(-1,0).$
For any $t\in\FF$ fix $\phi_t\in G$ such that 
$\phi_t(e_1)=te_1+e_2$ and $\phi_t(e_2)=e_1.$
The map
$$\FF\to(G\cdot (L',L))\cap(\Bc\times\{L\}),\quad 
t\mapsto(\phi_t(L),L)$$
is an isomorphism. 
The convolution product
$$\star\,:\,\CC_G[\Bc^n\times\Bc]\otimes\CC_G[\Bc\times\Bc]
\rightarrow\CC_G[\Bc^n\times\Bc]$$
is defined as
$$f\star g\, (L^n_\kb,L)=\sum_{L''\in\Bc}f(L^n_\kb,L'')\cdot g(L'',L).$$
Thus,
$$T_1(v_\jb)=q^{-1}\sum_{\kb\in\Pc_2}n(\jb,\kb)\,v_\kb,
\quad{\where}\quad
n(\jb,\kb)=\sharp\{t\in\FF\,|\,\phi_t^{-1}(L^n_\kb)\in B\cdot L^n_\jb\}.$$

\vskip1mm

\item{\bf .} Suppose first that $j_1=j_2$. In this case
$L_\jb^n$ is a sequence of lattices of the type $L(a,a)$ where $a\in\ZZ$. 
Thus, $L_\jb^n$ is fixed by $B$ and by $\phi_t$ for all $t$ and
$T_1(v_\jb)=q\,v_\jb.$

\vskip1mm

\item{\bf .} Suppose that $j_1<j_2$. Then $L_\jb^n$ is a sequence
of lattices of the type $L(a,b)$ with $a\geq b$ and the inequality is strict
for at least one lattice in the sequence. The only possibility to get
$L_\kb^n\in\phi_t(B\cdot L_\jb^n)$ for some $\kb\in\Pc_m$ and some $t\in\FF$
is that $t=0$ and, then, necessarily $\kb=(j_2,j_1)$.
Thus, $T_1(v_\jb)=q^{-1}\,v_{j_2j_1}.$

\vskip1mm

\item{\bf .} If $j_1>j_2$ the formula for $T_1(v_\jb)$ follows from 
the previous case and the relation $(T_1+q^{-1})(T_1-q)$. 

\vskip2mm

\noindent As for the $Q$ it is immediate that 
$\vartheta$ is the convolution product by the characteristic function
of the $G$-orbit of the pair $(L'',L)$, where $L''_i=L_{i+1}$
for all $i\in\ZZ$. Let us now compute the $\Uaf$-action. The integer
$m$ is no longer supposed to be 2. 
We consider the convolution product 
$$\star\,:\,\CC_G[\Bc^n\times\Bc^n]\otimes\CC_G[\Bc^n\times\Bc]
\rightarrow\CC_G[\Bc^n\times\Bc].$$
Given $\jb,\kb\in\Pc^0_m$ let
$\chi_{\jb,\kb}\in\CC_G[\Bc^n\times\Bc^n]$ be the characteristic function
of the $G$-orbit of the pair $(L^n_\jb,L^n_\kb)$. By definition the
map $\Uaf\rightarrow\CC_G[\Bc^n\times\Bc^n]$ send 
$\eb_i$, $\fb_i$ and $\kb_i$ respectively to
$$\matrix
{\ds\sum_{\jb\in\Pc_m^0}q^{-\sharp\jb^{-1}(i)}\,
\delta(\jb^{-1}(i+1)\neq\emptyset)\,\chi_{\jb^-_{s+1},\jb},}\hfill\cr\cr
{\ds\sum_{\jb\in\Pc_m^0}q^{-\sharp\jb^{-1}(i+1)}\,
\delta(\jb^{-1}(i)\neq\emptyset)\, \chi_{\jb^+_s,\jb},}\hfill\cr\cr
{\ds\sum_{\jb\in\Pc_m^0}q^{\sharp\jb^{-1}(i)-\sharp\jb^{-1}(i+1)}\,
\chi_{\jb,\jb},}\hfill
\endmatrix$$
where $s=\max\jb^{-1}(i)$ and $i=0,1,...,n-1$.
Suppose first that $\jb\in\Pc^0_m$. Then 
$$\Oc_\jb\cap(\Bc^n\times\{L\})=\{(L_\jb^n,L)\}.$$
If $\jb^{-1}(i+1)\neq\emptyset$ then
$$(G\cdot(L^n_{\jb^-_{s+1}},L^n_\jb))\cap(\Bc^n\times\{L_\jb^n\})=
\{L^n_{\jb^-_k}\,|\,k\in\jb^{-1}(i+1)\}\times\{L_\jb^n\},$$
and, thus,
$$\eb_i(v_\jb)=q^{-\sharp\jb^{-1}(i)}\,\sum_{k\in\jb^{-1}(i+1)}v_{\jb^-_k}.$$
The operator
$$v_\jb\mapsto q^{-\sharp\jb^{-1}(i)}\,
\sum_{k\in\jb^{-1}(i+1)}q^{2\sharp\{l\in\jb^{-1}(i)\,|\,l>k\}}v_{\jb^-_k},
\quad\forall\jb\in\Pc_m,$$
commutes to $\Haf$. When $\jb\in\Pc_m^0$ it is precisely the action
of $\eb_i$ written above. Thus it coincides with $\eb_i$. The proof
is similar for $\fb_i$ and $\kb_i$.
\qed

\vskip3mm

\noindent{\bf 7.}
The affine quantum group $\Uaf$ is known to admit
the structure of a Hopf algebra whose coproduct
$\Delta$ is such that for any $i=0,1,...,n-1$
$$\Delta(\eb_i)=\eb_i\otimes\kb_i+1\otimes\eb_i,\qquad
\Delta(\fb_i)=\fb_i\otimes 1+\kb_i^{-1}\otimes\fb_i,\qquad
\Delta(\kb_i)=\kb_i\otimes\kb_i.$$
As a consequence, $\Uaf$ acts 
on $\bigotimes^m\CC^n[\zeta^{\pm 1}]$
by iterating the action of $\eb_i,\fb_i,\kb_i$
on $\CC^n[\zeta^{\pm 1}]$ given in section 4.
Let us call this representation the {\it tensor representation} of $\Uaf$.
It is important to notice that the resulting representation of 
$\Uaf$ is isomorphic to the geometric one given in section 5.
The purpose of this section is to write explicitely such an isomorphism.
In particular it will follows that the formulas in section 5
do define a representation of $\Uaf$. 
Let us first prove the following technical result.
According to the proposition 1 the action of 
$\Haf$ on $\bigotimes^m\CC^n[\zeta^{\pm 1}]$ described in section 5
restricts to a representation
of the ring $\CC[X_1^{\pm 1},X_2^{\pm 1},...,X_m^{\pm 1}]$.

\vskip3mm

\noindent{\bf Lemma 7.} {\it The space $\bigotimes^m\CC^n[\zeta^{\pm 1}]$ is a 
free module over the ring $\CC[X_1^{\pm 1},X_2^{\pm 1},...,X_m^{\pm 1}]$, with
basis the monomials $v_\jb$ such that $\jb$ belongs to the set
$$\Pc_m^1=\{\jb\in\Pc_m\,|\,1\leq j_1, j_2,\cdots, j_m\leq n\}.$$
} 

\vskip1mm

\noindent{\it Proof.}
Let us consider $q$ as a formal variable.
For any $\jb\in\Pc_m$ we write 
$\jb=n\,\ujb+\ojb$, with $\ujb\in\Pc_m$ and $\ojb\in\Pc_m^1$.
Let $\Fb$ be the free $\CC[q^{\pm 1}]$-module with basis
$\Pc_m$. Let $p$ be the map
$$p\quad:\quad\Fb\to\bigoplus_{\jb\in\Pc_m}\CC[q^{\pm 1}]\,v_\jb,\qquad
\jb\mapsto X^{-\ujb}(v_\ojb).$$
The map $p$ is surjective. Namely for all $s,t,\jb$, there exists
a monomial $x$ in the $T^{\pm 1}_i$'s such that $\sigma_{s\,t}v_\jb=x(v_\jb)$.
Thus, for all $\jb$ there is a monomial $x$ in $Q^{\pm 1}$ and the 
$T_i^{\pm 1}$'s
such that $v_\ojb=x(v_\jb)$. 
The surjectivity follows since the $X^\ab T_w$'s
form a basis of $\Haf$.
As for the injectivity, the kernel of $p$ is free and vanishes when $q=1$. 
Thus we are done. 
\qed

\vskip3mm

\noindent{\bf Remark 7.1.} 
As a consequence if
$\Pc_m^0=\{\jb\in\Pc_m\,|\,1\leq j_1\leq j_2\leq\cdots\leq j_m\leq n\}$,
then $\bigotimes^m\CC^n[\zeta^{\pm 1}]= \sum_{\jb\in\Pc_m^0}\Haf\cdot v_\jb.$

\vskip3mm

\noindent Let $\Psi$ be the unique 
$\CC[X_1^{\pm 1},X_2^{\pm 1},...,X_m^{\pm 1}]$-linear automorphism
of $\bigotimes^m\CC^n[\zeta^{\pm 1}]$ such that 
$$\Psi(v_\jb)=q^{\sharp\{1\leq s<t\leq m\,|\,j_s<j_t\}}v_\jb,\qquad
\forall\jb\in\Pc_m^1,$$
and let $\dot\Phi$ be the linear isomorphism 
$$\dot\Phi\,:\quad
\bigotimes^m\CC^n[X^{\pm 1}]\,{\buildrel\sim\over\rightarrow}
\bigotimes^m\CC^n[\zeta^{\pm 1}],\qquad
X^\ab\otimes v_\jb\mapsto X^\ab\cdot\Psi(v_\jb),\qquad
\forall \jb\in\Pc^1_m.$$
The following result is stated without a proof in [GRV] and 
the analogue in the finite case is given in [GL].

\vskip3mm

\noindent{\bf Proposition 7.} {\it
For any $i=0,1,...,n-1$,
the operators $\dot\Phi^{-1}\circ\eb_i\circ\dot\Phi$, 
$\dot\Phi^{-1}\circ\fb_i\circ\dot\Phi$ and 
$\dot\Phi^{-1}\circ\kb_i\circ\dot\Phi$
on $\bigotimes^m\CC^n[X^{\pm 1}]$
do coincide with $\Delta^{m-1}(\eb_i)$, $\Delta^{m-1}(\fb_i)$,
and $\Delta^{m-1}(\kb_i)$.
Moreover for all $\jb\in\Pc_m^1$ and all 
$P\in\CC[X_1^{\pm 1}, X_2^{\pm 1},...,X_m^{\pm 1}]$,
$$\dot\Phi^{-1}\circ T_k\circ\dot\Phi (P\,v_\jb)=
\cases
qP^{s_k}\,v_\jb+
(q^{-1}-q){X_{k+1}(P-P^{s_k})\over X_k-X_{k+1}}\,v_\jb\quad
&\text{if $\quad j_k=j_{k+1},$}\cr\cr
P^{s_k}\,\sigma_kv_\jb+
(q^{-1}-q){X_{k+1}(P-P^{s_k})\over X_k-X_{k+1}}\,v_\jb\quad
&\text{if $\quad j_k<j_{k+1},$}\cr\cr
P^{s_k}\,\sigma_kv_\jb+
(q^{-1}-q){X_{k+1}P-X_kP^{s_k}\over X_k-X_{k+1}}\,v_\jb\quad
&\text{if $\quad j_k>j_{k+1},$}
\endcases$$
where $\sigma_k$ stands for $\sigma_{k,k+1}$ and $P^{s_k}$ is $P$ with
$X_k$ and $X_{k+1}$ exchanged.
}

\vskip3mm

\noindent{\it Proof.} Since $T_k$ commutes with any polynomial symmetric in 
the variables $X_k$ and $X_{k+1}$, we get 
$$\matrix
2T_k(P\,v_\jb)&=T_k((P+P^{s_k})+Q(X_k-X_{k+1}))\,v_\jb,\hfill\cr\cr
&=(P+P^{s_k})T_k\,v_\jb+QT_k(X_k-X_{k+1})\,v_\jb,\hfill
\endmatrix$$ 
where $Q$ is a symmetric polynomial in $X_k$ and $X_{k+1}.$ 
Now 
$$\matrix
T_kX_k=X_{k+1}T_k+(q^{-1}-q)X_{k+1},\hfill\cr\cr
-T_kX_{k+1}=(q^{-1}-q)X_{k+1}-X_kT_k,\hfill
\endmatrix
$$ 
from which we get 
$$T_k(P\,v_\jb)=P^{s_k}\,T_k\,v_\jb+(q^{-1}-q)\,{X_{k+1}(P-P^{s_k})\over
X_k-X_{k+1}}\,v_\jb.$$ 
Thus
$$\dot\Phi^{-1}\circ T_k\circ\dot\Phi (P\,v_\jb)= 
q^{\sharp\{1\leq s<t\leq m|j_s<j_t\}}P^{s_k}\Psi^{-1}
(T_k\,v_\jb)+(q^{-1}-q){X_{k+1}(P-P^{s_k})\over X_k-X_{k+1}}\, v_\jb$$ 
and the formulas follow from an easy case by case computation,
according to the value of $T_kv_\jb$ as in section 5.
As for $\Uaf$, let consider the case of $\eb_i$, $i=0,1,...,n-1$,
since the other cases are quite similar. Then

$$\matrix
{\ds \dot\Phi^{-1}\circ\eb_i\circ\dot\Phi (P\,v_\jb)=}\hfill\cr\cr
{\ds =q^{\sharp\{s<t\,|\,j_s<j_t\}-\sharp\jb^{-1}(i)}
\sum_{k\in\jb^{-1}(i+1)}q^{2\sharp\{l\in\jb^{-1}(i)\,|\,l>k\}}
P\,\Psi^{-1}(v_{j_k^-})},\hfill\cr\cr
{\ds =\sum_{k\in\jb^{-1}(i+1)}q^{2\sharp\{l\in\jb^{-1}(i)\,|\,l>k\}
-\sharp\jb^{-1}(i)
-\sharp\{l\in\jb^{-1}(i+1)\,|\,l>k\}
+\sharp\{l\in\jb^{-1}(i)\,|\,l<k\}}
P\,v_{j_k^-}},\hfill\cr\cr
{\ds=\sum_{k\in\jb^{-1}(i+1)}q^{\sharp\{l\in\jb^{-1}(i)\,|\,l>k\}
-\sharp\{l\in\jb^{-1}(i+1)\,|\,l>k\}}
P\,v_{j_k^-}},\hfill\cr\cr
{\ds=\sum_{k=1}^m(1^{\otimes k-1}\otimes\eb_i\otimes\kb_i^{\otimes
m-k})\,v_\jb},\hfill\cr\cr
{\ds=(\Delta^{m-1}\eb_i)\,v_\jb.}\hfill
\endmatrix$$
\qed

\vskip3mm

\noindent{\bf Remark 7.2.} 
The action of the $T_i$'s on $\bigotimes^m\CC^n[X^{\pm 1}]$
given above and the product by the $X_i$'s determine a representation of 
$\Haf$ on the tensor module, introduced for the first time in [GRV]. 

\vskip3mm

\noindent{\bf 8.} Let now consider the K-theoretic analogue of the previous
construction in the same way as in [GRV], [GKV]. 
Fix another set of formal variables $z_1^{\pm 1},...,z_m^{\pm 1}$.
The purpose of this section is to explain how the
actions of $\Haf$ and $\Uaf$ on $\bigotimes^m\CC^n[\zeta^{\pm 1}]$
defined in section 5 can be induced to commuting representations
of $\Htor$ and $\Utor$ on the space
$\Vb_m=(\CC^n)^{\otimes m}
[\zeta_1^{\pm 1},...,\zeta_m^{\pm 1},z_1^{\pm 1},...,z_m^{\pm 1}]$.
The formulas for the induced action of $\Htor$ may be taken from [C2] 
for instance. The generators $T_i,Y_i,Q$ act as follows :
for any $P\in\Rb_m$ and $v\in\bigotimes^m\CC^n[\zeta^{\pm 1}]$,
$$\matrix
T_i(v\cdot P)=(\tau_{i,i+1}(v)-qv)\cdot s_{i,i+1}(P)+v\cdot t_{i,i+1}(P),
\hfill\cr\cr
Y_i(v\cdot P)=v\cdot Pz_i^{-1},\hfill\cr\cr
Q(v\cdot P)=\vartheta(v)\cdot D_1s_{1,m}s_{1,m-1}\cdots s_{1,2}(P),\hfill
\endmatrix$$
and the central element $\xb$ goes to a fixed $p\in\CC^\times$.
The operators $t_{i,i+1},s_{i,i+1}$ are defined in section 2 and 
$\tau_{i,i+1},\vartheta$ are defined in section 5.
As for the action of $\Utor$ on $\Vb_m$
let  $\eb_i,\fb_i,\kb_i\in\End_{\Htor}(\Vb_m)$  
($i=0,1,2,...,n-1$) be as in section 5 and
consider an additional triple of operators 
$\eb_n,\,\fb_n,\,\kb_n\in\End_{\Htor}(\Vb_m)$  
such that for all $\jb\in\Pc_m^0$,
$$\matrix
{\ds\eb_n(v_\jb)=q^{\sharp\jb^{-1}(n)-1}p^{1/n}\sum_{k\in\jb^{-1}(1)}
q^{2k-1-m}Y_k^{-1}\cdot v_{\jb_k^-}\zeta^{-1}_k,}\hfill\cr\cr
{\ds\fb_n(v_\jb)=q^{\sharp\jb^{-1}(1)+1}p^{-1/n}\sum_{k\in\jb^{-1}(n)}
q^{m-2k+1}Y_k\cdot v_{\jb_k^+}\zeta_k,}\hfill\cr\cr
{\ds\kb_n(v_\jb)=q^{\sharp\jb^{-1}(n)-\sharp\jb^{-1}(1)}v_\jb,}\hfill
\endmatrix$$
where $p^{1/n}$ is a fixed $n$-th root of $p.$

\vskip3mm

\noindent{\bf Proposition 8.} {\it
The operators $\eb_i,\,\fb_i,\,\kb_i$, $i=0,1,2,...,n-1$, 
(resp. $i=1,2,...,n$) 
define an action of $\Uaf$ on $\Vb_m$ commuting with $\Htor$. 
Moreover these two actions of $\Uaf$ can be glued in a representation of
$\Utor$ commuting to $\Htor$ if $d=q^{-1}p^{1/n}$.}

\vskip3mm

\noindent 
Once again we do not prove this proposition here since it follows
immediately from the proof of proposition 9.

\vskip3mm

\noindent{\bf 9.} 
Let us recall that
$\Vb_m=(\CC^n)^{\otimes m}
[\zeta_1^{\pm 1},...,\zeta_m^{\pm 1},z_1^{\pm 1},...,z_m^{\pm 1}]$.
In section 8 we have defined an action of $\Htor$ and $\Utor$ on $\Vb_m$.
The Schur duality is an equivalence of categories between finite dimensional
representations of the symmetric group and of the linear group. It has been
generalized to quantum groups by Jimbo, Drinfeld and Cherednik. In [VV]
we proved a similar duality between $\Htor$ and $\Utor$.
Since the $\Utor$-action on $\Vb_m$ commutes to $\Htor$, to any
right ideal $\Jb\subset\Htor$ we can associate a representation
of the quantized toroidal algebra on the quotient
$(\Jb\cdot\Vb_m)\backslash\Vb_m$. 
Let consider the right $\Htor$-module
$\Jb\backslash\Htor$. The purpose of this section
is to prove that the Schur dual of $\Jb\backslash\Htor$, as defined in [VV], 
is isomorphic to $(\Jb\cdot\Vb_m)\backslash\Vb_m$. 
By definition, the underlying vector space of
the Schur dual of $\Jb\backslash\Htor$ is
$$(\Jb\backslash\Htor)\otimes_{\Hb_m}(\CC^n)^{\otimes m}.$$
In order to describe the action of $\Utor$ on this space set
$$\eb_\theta(v_j)=\delta(j=n)v_1,\qquad
\fb_\theta(v_j)=\delta(j=1)v_n,\qquad
\kb_\theta(v_j)=q^{\delta(j=1)-\delta(j=n)}v_j,$$
and 
$\fb_{\theta,l}=1^{\otimes
l-1}\otimes\fb_\theta\otimes(\kb_\theta^{-1})^{\otimes m-l}$,
$\,\eb_{\theta,l}=\kb_\theta^{\otimes
l-1}\otimes\eb_\theta\otimes 1^{\otimes m-l}$,
for all $j=1,2,...,n$ and $l=1,2,...,m$. Then for all $\jb\in\Pc_m^1$ and
$x\in\Jb\backslash\Htor$,
$$\matrix
{\ds\eb_0(x\otimes v_\jb)=
\sum_{l=1}^mx\,X_l\otimes\fb_{\theta,l}(v_\jb)},\qquad\hfill &
{\ds \fb_0(x\otimes v_\jb)=
\sum_{l=1}^mx\,X_l^{-1}\otimes\eb_{\theta,l}(v_\jb),}
\hfill\cr\cr
{\ds\eb_n(x\otimes v_\jb)=
q^{-1}p^{1/n}\sum_{l=1}^mx\,Y^{-1}_l\otimes \fb_{\theta,l}(v_\jb)},
\qquad\hfill &
{\ds\fb_n(x\otimes v_\jb)=
qp^{-1/n}\sum_{l=1}^mx\,Y_l\otimes\eb_{\theta,l}(v_\jb),}\hfill
\endmatrix$$
and, if $i=1,2,...,n-1$,
$$\eb_i=1\otimes\Delta^{m-1}\eb_i,\qquad
\fb_i=1\otimes\Delta^{m-1}\fb_i$$

\noindent (see [VV] and [CP] for more details). In the particular case
$\Jb=\{0\}$, the Schur dual of the right regular module $\Htor$
is endowed with an action of $\Htor$ by left translations which commutes
to $\Utor$.

\vskip3mm

\noindent{\bf Proposition 9.} {\it Fix $d=q^{-1}p^{1/n}$.
The $\Utor\times\Htor$-module $\Vb_m$ is isomorphic to 
the Schur dual of the right regular representation of $\Htor$.
As a consequence the quotient $(\Jb\cdot\Vb_m)\backslash\Vb_m$ 
is isomorphic to the Schur dual of $\Jb\backslash\Htor$ for any right ideal 
$\Jb\subset\Htor$.} 

\vskip3mm

\noindent{\it Proof.} 
The proposition 7 implies that the bijection 
$${\ts
\dot\Phi\,:\quad
\Haf\otimes_{\Hb_m}(\CC^n)^{\otimes m}\,{\buildrel\sim\over\rightarrow}
\bigotimes^m\CC^n[X^{\pm 1}]\,{\buildrel\sim\over\rightarrow}
\bigotimes^m\CC^n[\zeta^{\pm 1}],}$$
intertwines the tensor representation of $\Uaf$ and $\Haf$ given
in section 7, and the convolution
action of $\Uaf$ and $\Haf$ described in section 5. 
This isomorphism extends to a bijection $\ddot\Phi$ 
$$\ddot\Phi\,:\quad\Htor\otimes_{\Hb_m}(\CC^n)^{\otimes m}
\,{\buildrel\sim\over\rightarrow}\Vb_m,\quad
Y^\bb X^\ab \otimes v\mapsto Y^\bb X^\ab\cdot\Psi(v),\quad
\forall v\in(\CC^n)^{\otimes m},$$
where $\Htor$ acts on $\Vb_m$ as in section 8.
This morphism is well defined because the monomials $Y^\bb X^\ab$ form
a basis of $\Htor$ as a right $\Hb_m$-module.
Now,

\vskip2mm

\itemitem{\bf .}
$\Vb_m$ is the $\Htor$-module
induced from the $\Haf$-module $\bigotimes^m\CC^n[\zeta^{\pm 1}]$
(see section 8),

\vskip2mm

\itemitem{\bf .}
$\Htor\otimes_{\Hb_m}(\CC^n)^{\otimes m}$ is the $\Htor$-module
induced from the $\Haf$-module
$\Haf\otimes_{\Hb_m}(\CC^n)^{\otimes m}$,

\vskip2mm

\itemitem{\bf .}
$\dot\Phi$ is an isomorphism of $\Haf$-modules.

\vskip2mm

\noindent Thus $\ddot\Phi$ is an isomorphism of $\Htor$-modules.
As for the $\Utor$-action let us do a direct computation.
Since $\ddot\Phi$ commutes with the left action of $\Htor$ it
suffices to prove that for all $\jb\in\Pc^0_m$
and all $i=0,1,...,n$,
$$\ddot\Phi\circ\eb_i(1\otimes v_\jb)=\eb_i\circ\ddot\Phi(1\otimes v_\jb)
\quad\and
\quad\ddot\Phi\circ\fb_i(1\otimes v_\jb)=\fb_i\circ\ddot\Phi(1\otimes v_\jb).$$
If $i\neq n$ this has already been proved in the proposition 7.
As for the remaining cases, we get
$$\matrix
{\ds\ddot\Phi\circ\eb_n(1\otimes v_\jb)}
&{\ds =q^{-1}p^{1/n}\sum_{l=1}^m\,
\ddot\Phi(Y^{-1}_l\otimes \fb_{\theta,l}(v_\jb))},\hfill\cr\cr
&{\ds =q^{-1}p^{1/n}\sum_{l=1}^m\,
Y^{-1}_l\cdot\Psi(\fb_{\theta,l}(v_\jb))},\hfill\cr\cr
&{\ds =q^{-1}p^{1/n}\sum_{l\in\jb^{-1}(1)}\,
q^{\sharp\jb^{-1}(n)-\sharp\jb^{-1}(1)+l}
Y^{-1}_l\cdot\Psi(v_{\jb_l^-}\zeta_l^{-1})},\hfill\cr\cr
&{\ds =q^{\sharp\jb^{-1}(n)-1+\sharp\{s<t\,|\,j_s<j_t\}}p^{1/n}
\sum_{l\in\jb^{-1}(1)}\,
q^{2l-m-1}Y^{-1}_l\cdot v_{\jb_l^-}\zeta_l^{-1}},\hfill\cr\cr
&{\ds=\eb_n\circ\ddot\Phi(1\otimes v_\jb).}\hfill
\endmatrix
$$
The computation for $\fb_n$ is similar.\qed

\vskip3mm

\noindent{\bf 10.} 
In the remaining three sections we fix $d=q^{-1}p^{1/n}$.
Let consider the right ideal 
$\gamma\omega(\Ib_m)\subset\Htor$, where $\Ib_m$ is the ideal defined 
in section 2.
Let $\bigwedge^m$ be the quotient 
$(\gamma\omega(\Ib_m)\cdot\Vb_m)\backslash\Vb_m$. 
Recall that $\Haf$ acts on the tensor module
$$\bigotimes^m\CC^n[X^{\pm 1}]=
(\CC^n)^{\otimes m}[X_1^{\pm 1},X_2^{\pm1},...,X_m^{\pm 1}]$$
as in the proposition 7. Then, set (see [KMS])
$$\Omega=\sum_{i=1}^{m-1}\Ker(T_i-q)
\subset\bigotimes^m\CC^n[X^{\pm 1}].$$

\vskip3mm

\noindent{\bf Lemma 10.} {\it
The space $\bigwedge^m$ is isomorphic to the quotient
$\bigotimes^m\CC^n[X^{\pm 1}]/\Omega.$
}

\vskip3mm

\noindent{\it Proof.}
By definition, $\bigwedge^m$ is isomorphic to 
$$(\gamma\omega(\Ib_m)\backslash\Htor)\otimes_{\Hb_m}(\CC^n)^{\otimes m}.$$
The left ideal $\Ib_m$ is generated by the $T_i-q$'s and by $Q-1$. Thus
$\gamma\omega(\Ib_m)$ is the right ideal of $\Htor$ generated by the
$T_i+q^{-1}$'s and the $Y_i-q^{2i-2}$'s. Since the monomials 
$Y^\bb\,X^\ab\,T_w$, 
$k\in\ZZ$, $\ab,\bb\in\ZZ^m$, $w\in{\frak{S}}_m$, form a basis of $\Htor$,
the map
$$(\gamma\omega(\Ib_m)+Y^\bb\,X^\ab\,T_w)\otimes v
\mapsto q^{2\sum_i(i-1)b_i}X^{\ab}\,(T_wv),$$
for all $k\in\ZZ$, $\ab\in\ZZ^m$ and $v\in(\CC^n)^{\otimes m}$,
is an isomorphism from the space of $q$-wedges to the quotient of
$\bigotimes^m\CC^n[X^{\pm 1}]$
by the right ideal generated by the $T_i+q^{-1}$'s. 
We are thus reduced to prove that
$$\sum_{i=1}^{m-1}\Ker(T_i-q)=\sum_{i=1}^{m-1}\Im(T_i+q^{-1})$$
in $\bigotimes^m\CC^n[X^{\pm 1}]$.
One inclusion follows from the relation $(T_i-q)(T_i+q^{-1})=0$.
The equality can be proved by using the formulas in section 5 since
$$(\tau_{1,2}-q)(v_{ij})=\cases
0&\text{if $\quad i=j,$}\cr
q^{-1}v_{ji}-qv_{ij}&\text{if $\quad i<j,$}\cr
qv_{ji}-q^{-1}v_{ij}&\text{if $\quad i>j,$}\cr
\endcases$$
and
$$(\tau_{1,2}+q^{-1})(v_{ij})=\cases
(q+q^{-1})v_{ij}&\text{if $\quad i=j,$}\cr
q^{-1}(v_{ji}+v_{ij})&\text{if $\quad i<j,$}\cr
q(v_{ji}+v_{ij})&\text{if $\quad i>j,$}\cr
\endcases$$
for all $i,j\in\ZZ$ (where, as in section 4, $v_{i+nk}=v_i\,\zeta^{-k}$ for all
$k\in\ZZ$ and $i=1,2,...,n$).
\qed

\vskip3mm

\noindent 
According to [KMS], elements 
of $\bigwedge^m$ are called {\it $q$-wedges}. 
From now on it will be more convenient to view the $q$-wedges space as the
quotient $\bigotimes^m\CC^n[X^{\pm 1}]/\Omega.$ 
Let $\wedge\,:\,\bigotimes^m\CC^n[X^{\pm 1}]\,\to\,\bigwedge^m$
be the projection. For all $\jb\in\Pc_m$ we denote indifferently by
$$\wedge v_\jb\quad{\or}\quad 
v_{j_1}\wedge v_{j_2}\wedge\cdots\wedge v_{j_m}$$
the projection of the monomial
$v_\jb=v_{j_1\,j_2\cdots j_m}$ in the $q$-wedges space.
From section 9 the algebra $\Utor$ acts on $\bigwedge^m$.
Let $\hat x_l$ be the image of $\omega(X_l)$ in $\End(\Rb_m)$,
where $X_l$ acts on $\Rb_m$ by the operator $x_l$ defined in section 2.

\vskip3mm

\noindent{\bf Theorem 10.} {\it The action of the generators of $\Utor$ on
$\bigwedge^m$ is given as follows : $\forall\jb\in\Pc_m^1,$
$\forall P(X)\in\CC[X_1^{\pm 1},X_2^{\pm 1},...,X_m^{\pm 1}],$
$\forall i=1,2,...,n-1,$
$$\matrix
{\ds\eb_i\cdot\wedge P(X^{-1})\,v_\jb=
\wedge(P(X^{-1})\,\Delta^{m-1}(\eb_i)\cdot v_\jb),}\hfill\cr\cr
{\ds\fb_i\cdot\wedge P(X^{-1})\,v_\jb=
\wedge(P(X^{-1})\,\Delta^{m-1}(\fb_i)\cdot v_\jb),}\hfill\cr\cr
{\ds\kb_i\cdot\wedge P(X^{-1})\,v_\jb=
\wedge(P(X^{-1})\,\Delta^{m-1}(\kb_i)\cdot v_\jb),}\hfill\cr\cr
{\ds\eb_0\cdot\wedge P(X^{-1})\,v_\jb=
\sum_{l=1}^m\wedge(X_lP(X^{-1})\,\fb_{\theta,l}\cdot v_\jb),}
\hfill\cr\cr
{\ds\fb_0\cdot\wedge P(X^{-1})\,v_\jb=
\sum_{l=1}^m\wedge(X_l^{-1}P(X^{-1})\,\eb_{\theta,l}\cdot v_\jb),}
\hfill\cr\cr
{\ds\kb_0\cdot\wedge P(X^{-1})\,v_\jb=
\wedge(P(X^{-1})\,\kb_{\theta}^{-1\,\otimes l}\cdot v_\jb),}
\hfill\cr\cr
{\ds\eb_n\cdot\wedge P(X^{-1})\,v_\jb=q^{-1}p^{1/n}
\sum_{l=1}^m\wedge(\hat x_l(P)(X^{-1})\,\fb_{\theta,l}\cdot v_\jb),
}\hfill\cr\cr
{\ds\fb_n\cdot\wedge P(X^{-1})\,v_\jb=qp^{-1/n}
\sum_{l=1}^m\wedge(\hat x_l^{-1}(P)(X^{-1})\,\eb_{\theta,l}\cdot v_\jb),
}\hfill\cr\cr
{\ds\kb_n\cdot\wedge P(X^{-1})\,v_\jb=
\wedge(P(X^{-1})\,\kb_{\theta}^{-1\,\otimes l}\cdot v_\jb).}
\hfill\cr\cr
\endmatrix$$
}

\vskip3mm

\noindent{\it Proof.}
According to the isomorphism in the proposition above and the description 
of the Schur dual
recalled in section 9, we must prove that
$P(X^{-1})Y_l^{-1}-\hat x_l(P)(X^{-1})\in\gamma\omega(\Ib_m)$ 
for all $l$ and $P$.
By definition of $\hat x_l$ we have
$\omega(X_l)P(Y^{-1})-\hat x_l(P)(Y^{-1})\in\Ib_m$.
Thus
$$\omega\gamma(P(X^{-1})Y_l^{-1}-\hat x_l(P)(X^{-1}))=
\omega(X_l)P(Y^{-1})-\hat x_l(P)(Y^{-1})\in\Ib_m.$$
\qed

\vskip3mm

\noindent{\bf Remark 10.}  A direct computation gives
$$\hat x_l=q^{m-1}\,t_{_{l-1,l}}^{-1}\,s_{_{l-1,l}}\,t_{_{l-2,l}}^{-1}\,
s_{_{l-2,l}}\cdots t_{_{1,l}}^{-1}\,s_{_{1,l}}\,D_l\,s_{_{l,m}}\,t_{_{l,m}}
\cdots s_{_{l,l+1}}\,t_{_{l,l+1}}.$$
In particular $\hat x_l=D_l$ if $q$ is one.

\vskip3mm

\noindent
Let us notice that it is possible to write explicit formulas for 
the action of all the Drinfeld generators on the tensor module. 
We will use these formulas in the proof of theorem 12.
Given $1\leq i < j \leq m$, set
$$T_{i,j}=T_iT_{i+1}\cdots T_j\quad{\and}\quad
T_{j,i}=T_jT_{j-1}\cdots T_i.$$

\vskip3mm

\noindent{\bf Proposition 10. [VV; Theorem 3.3]} 
{\it For any $\jb\in\Pc_m^0$ and $i=1,2,...,n-1$, set
$s=\max\jb^{-1}(i)$. Then,

$$
\fb_i(w)\cdot\wedge P(X^{-1})\,v_\jb=
q^{1-\sharp\jb^{-1}(i)}
\sum_{k\in\jb^{-1}(i)}q^{^{s-k}}\wedge P(X^{-1})T_{_{k,s-1}}
\epsilon(q^np^{i/n}wY_{_s}) v_{\jb_s^+},$$
if  $\,\jb^{-1}(i)\not=\emptyset$, $0$ else, and  
$$\eb_i(w)\cdot\wedge P(X^{-1})\,v_\jb=
q^{1-\sharp\jb^{-1}(i+1)}
\sum_{k\in\jb^{-1}(i+1)}q^{^{k-s-1}}\wedge P(X^{-1})T_{_{k-1,s+1}}
\epsilon(q^np^{i/n}wY_{_{s+1}}) v_{\jb_{s+1}^-},$$
if  $\,\jb^{-1}(i+1)\not=\emptyset$, $0$ else.
}
\vskip1mm

\noindent The shift in the powers of $q$ in [VV] is due to a different 
normalization of the $T_i$'s.

\vskip3mm

\noindent{\bf 11.} 
According to [KMS] set $u_\jb=X^{-\ujb}v_{\ojb}$ for all $\jb\in\Pc_m$.
The Fock space, $\bigwedge^{\infty/2}$, is the linear span of semi-infinite 
monomials 
$$\wedge u_\jb=u_{j_1}\wedge u_{j_2}\wedge u_{j_3}\wedge u_{j_4}\wedge\cdots$$
where $j_1>j_2>j_3>j_4>...$ and $j_{k+1}=j_k-1$ for $k>>1$. 
It splits as a direct sum of an infinite number of {\it sectors},
$\bigwedge^{\infty/2}_{(e)}$, parametrized by an integer $e\in\ZZ$. 
Let denote by $\eb=(e_1,e_2,e_3...)$ the infinite sequence such 
that $e_k=e-k+1$ for any $k\in\ZZ$. Then,
$\bigwedge^{\infty/2}_{(e)}$ is the linear span of semi-infinite 
monomials $\wedge u_\jb$
such that $j_1>j_2>j_3>j_4>...$ and $j_k=e_k$ for $k>>1$. 
In particular the sector $\bigwedge^{\infty/2}_{(e)}$ contains the element 
$$|e\rangle=u_{e_1}\wedge u_{e_2}\wedge u_{e_3}\wedge u_{e_4}\wedge\cdots$$
called the vacuum vector. 
The space $\bigwedge^{\infty/2}_{(e)}$ is $\NN$-graded. 
Recall that the height of an infinite sequence $\jb=(j_1,j_2,j_3,...)$
with only a finite number of non-zero terms is defined as $|\jb|=\sum_kj_k$.
Then, if $j_1>j_2>j_3>j_4>...$ and $j_k=e_k$ for $k>>1$, set
$$\deg(\wedge u_\jb)=|\ujb-\ueb|.$$
The degree of the vacuum vector is zero. 
We denote by $\bigwedge^{\infty/2,k}_{(e)}\subset\bigwedge^{\infty/2}_{(e)}$ 
the homogeneous component of degree $k$.

\vskip3mm

\noindent{\bf 12.}
Fix $e\in\ZZ$, $p\in\CC^\times$ and a generic $q\in\CC^\times$.
As before, $m$ is a non-negative integer.
We want to construct a representation of $\Utor$ on the Fock space as a limit 
when $m\to\infty$ of
the representation of $\Utor$ on $\bigwedge^m$.
It is proved in [KMS] that the {\it normally ordered} $q$-wedges, i.e.
the $\wedge u_\jb$'s such that
$$\jb\in\Pc_m^{no}=\{\jb\in\Pc_m\,|\,j_1>j_2>\cdots>j_m\},$$
form a basis of $\bigwedge^m$. 
Put $\eb^{m}=(e_1,e_2,...,e_m)$. 
Let $\bigwedge^m_{(e)}\subset\bigwedge^m$ be the linear span
of the normally ordered $q$-wedges $\wedge u_\jb$ such that 
$\ujb\geq\ueb^{m},$ i.e. $\uj_k\geq\ue_k$ for all $k$.
The space $\bigwedge^m_{(e)}$ admits a grading similar to the Fock space 
grading : 
$$\deg(\wedge u_\jb)=|\ujb-\ueb^{m}|,
\qquad\forall\jb\in\Pc_m^{no}\quad{\s.\t.}\quad\ujb\geq\ueb^{m}.$$
Let $\bigwedge^{m,k}_{(e)}\subset\bigwedge^{m}_{(e)}$ be the component of 
degree $k$ and set
$$\Pc_m^{no,k}=\{\jb\in\Pc_m^{no}\,|\,|\ujb-\ueb^m|=k\}.$$
Let
$\pi^{m,k}\,:\,\bigwedge^{m+n,k}_{(e)}\longrightarrow\bigwedge^{m,k}_{(e)}$
be the projections
$$\wedge u_\jb\mapsto\cases
\wedge u_{j_1\,j_2\cdots j_m}\quad&
\text{if $\quad\uj_k=\ue_l\quad\forall k,l\in[m+1,m+n]$,}\cr
0&\text{else,}\endcases$$
for all $\jb\in\Pc^{no,k}_{n+m}$.

\vskip3mm

\noindent {\bf Proposition 12. [TU]} {\it For any $m,k\in\NN$, 

\vskip1mm

\noindent {\rm (12.1)}\, $\bigwedge^{m,k}_{(e)}\subset\bigwedge^m$ is 
$\Uafv$-stable. 

\vskip1mm

\noindent {\rm (12.2)}\, If $\om=\oe$ then
$\pi^{m,k}\in\Hom_{\Utor}(\bigwedge^{m+n,k}_{(e)},\bigwedge^{m,k}_{(e)}).$
If moreover $\um\geq k$ then $\pi^{m,k}$ is invertible.

\vskip1mm 

\noindent {\rm (12.3)}\,
Suppose that $\om=\oe$ and $\um\geq k$. Given $s\in\NN$ and 
$\jb=(\ib,\lb)\in\Pc^{no,k}_{m+s}$ with $\ib\in\Pc^{no}_m$ and $\lb\in
\Pc^{no}_s,$
$$
\left\{\matrix
\ui_m>\ul_{m+1}=\cdots =\ul_{m+s}\hfill\cr
1\leq r\leq m\hfill
\endmatrix\right.
\quad\Longrightarrow\quad
\wedge(X^{-\ujb}Y_r^{\pm 1}v_\ojb)=
\wedge(X^{-\uib}Y_r^{\pm 1}X^{-\ulb}v_\ojb).$$
}

\vskip1mm

\noindent{\it Proof.} The statement
(12.1) follows from [TU; Proposition 4 and (4.28)] and (12.2) follows from 
[TU; Proposition 5 and 6]. The formula (12.3) follows from 
[TU; (4.40)], from $\um\geq k$ and from the fact that both terms in the
equality have degree $k$.
The shift in the powers of $q$ in [TU; (4.40)] is due to a different 
normalization of the $Y_i$'s. 
\qed

\vskip3mm

\noindent As a consequence, if $\om=\oe$ and $\um\geq k$, the map
$${\ts\bigwedge^{m,k}_{(e)}\longrightarrow\bigwedge^{\infty/2,k}_{(e)},\qquad
v\mapsto v\wedge |e_{1+m}\rangle},$$
is a linear  isomorphism and $\bigwedge^{\infty/2,k}_{(e)}$ inherits a 
$\Uafv$-action from $\bigwedge^{m,k}_{(e)}$
which is independent of the choice of such an $m$.
For any $v\in\bigwedge^{m,k}_{(e)}$ and any $i=0,1,...,n$ set
$$\matrix
\eb_i(v\wedge|e_{1+m}\rangle)=\eb_i(v)\wedge\kb_i\,|e_{1+m}\rangle+
v\wedge\eb_i\,|e_{1+m}\rangle,\hfill\cr\cr
\fb_i(v\wedge|e_{1+m}\rangle)=\fb_i(v)\wedge|e_{1+m}\rangle+
\kb_i^{-1}(v)\wedge\fb_i\,|e_{1+m}\rangle,\hfill\cr\cr
\kb_i(v\wedge|e_{1+m}\rangle)=\kb_i(v)\wedge\kb_i\,|e_{1+m}\rangle,\hfill
\endmatrix$$
where
$$\eb_i\,|e_{1+m}\rangle=0,\quad
\fb_i\,|e_{1+m}\rangle=\delta(i=0)u_{e_{m}}\wedge|e_{2+m}\rangle,\quad
\kb_i\,|e_{1+m}\rangle=q^{\delta(i=0)}|e_{1+m}\rangle.
$$
Thus $\eb_i,\fb_i,\kb_i$, $i=1,2,...,n$, are precisely the generators of 
the action of $\Uafv$ on $\bigwedge^{\infty/2,k}_{(e)}$.

\vskip3mm

\noindent{\bf Theorem 12.} {\it 
The formulas above define a representation
of $\Utor$ on each sector of the Fock space.}

\vskip3mm

\noindent{\it Proof.} Let define $\phi_\infty$ as the linear automorphism
of the Fock space $\wedge u_\jb\mapsto\wedge u_{1+\jb}$,
for all normally ordered $q$-wedges $\wedge u_\jb\in\bigwedge^{\infty/2}$. 
Let us prove that
$$\matrix
\eb_i(w)=\phi_\infty^{-1}\circ \eb_{i+1}(p^{1/n}w)\circ\phi_\infty,\qquad\hfill
&\fb_i(w)=\phi_\infty^{-1}\circ\fb_{i+1}(p^{1/n}w)\circ\phi_\infty,\hfill\cr\cr
\kb^\pm_i(w)=\phi_\infty^{-1}\circ\kb^\pm_{i+1}(p^{1/n}w)\circ\phi_\infty,\hfill
&\forall i=1,2,..,n-2,\hfill
\endmatrix
\leqno(12.4)$$
$$\matrix
\eb_{n-1}(w)=\phi_\infty^{-2}\circ\eb_1(p^{2/n}w)\circ\phi_\infty^2,\qquad\hfill
&\fb_{n-1}(w)=\phi_\infty^{-2}\circ\fb_1(p^{2/n}w)\circ\phi_\infty^2,
\hfill\cr\cr
\kb^\pm_{n-1}(w)=
\phi_\infty^{-2}\circ\kb^\pm_1(p^{2/n}w)\circ\phi_\infty^2.\hfill&
\endmatrix
\leqno(12.5)$$
The algebra $\Uafv$ acts on each sector. Moreover, the map
$$\eb_{i,k}\mapsto a^k\eb_{i+1,k},\quad
\fb_{i,k}\mapsto a^k\fb_{i+1,k},\quad
\hb_{i,k}\mapsto a^k\hb_{i+1,k},$$
(where the index $n+1$ stands for 0)
extends to an automorphism of $\Utor$ for any $a\in\CC^\times$. 
As a consequence of (12.4) and (12.5), setting 
$$\eb_0(w)=\phi_\infty^{-1}\circ \eb_1(p^{1/n}w)\circ\phi_\infty,\qquad
\fb_0(w)=\phi_\infty^{-1}\circ\fb_1(p^{1/n}w)\circ\phi_\infty,\qquad
\kb^\pm_0(w)=\phi_\infty^{-1}\circ\kb^\pm_1(p^{1/n}w)\circ\phi_\infty,$$
we will get an action of $\Utor$ on the Fock space.
Since the operators $\eb_0,\fb_0$ and $\kb^\pm_0$ coincide with the 
degree zero Fourier components of the formal series
$\eb_0(w),\fb_0(w)$ and $\kb^\pm_0(w)$
respectively, we will be done. 
For any $m\in\NN$, consider the map
$${\ts\phi_m\quad:\quad\bigwedge^m_{(e)}\to\bigwedge^{m+1}_{(e+1)},\qquad
\wedge u_{\jb}\,\mapsto \wedge u_{1+\jb}\wedge u_{e_m},
\quad\forall\jb\in\Pc_{m}^{no}.}$$
Set $h=\oe+(n+1)k$ and
$\bigwedge^{\infty/2,\leq h}_{(e)}=
\bigoplus_{i\leq h}\bigwedge^{\infty/2,i}_{(e)}.$
Let us observe that
$$(12.2)\quad\Longrightarrow\quad\deg\phi_\infty(u)\leq h,\quad
\forall u\in{\ts\bigwedge^{\infty/2,k}_{(e)}}.$$
The maps $\phi_\infty$ and $\phi_m$ do not preserve the degree.
However, if $\om=\oe$ and $\um\geq h$, we get the following 
commutative diagram 
$$\matrix
&\bigwedge^{\infty/2,k}_{(e)}&{\buildrel\sim\over\longrightarrow}
&\bigwedge^{m,k}_{(e)}&\cr
{\ss\phi_\infty}\hskip-1cm&\downarrow&&\downarrow
&\hskip-1cm{\ss\phi_{m}}\cr
&\bigwedge^{\infty/2,\leq h}_{(e+1)}&{\buildrel\sim\over\longrightarrow}&
\bigwedge^{m+1,\leq h}_{(e+1)}&,
\endmatrix$$
where the horizontal arrows are the projections on the first $m$ (resp. $m+1$)
components.
Since the horizontal maps commute to $\Uafv$,
it suffices to prove (12.4), (12.5) on
$\bigwedge^{m,k}_{(e)}$ with respect to $\phi_m$.
As a consequence of the theorem 10 and [VV, proposition 3.4], the map
$${\ts\bigwedge^m\longrightarrow\bigwedge^m,\quad
\wedge u_\jb\mapsto\wedge u_{1+\jb},\quad
\forall \jb\in\Pc^{no}_m,}$$
intertwines $\eb_i(w),\fb_i(w),\kb_i^\pm(w)$ and
$\eb_{i+1}(w),\fb_{i+1}(w),\kb_{i+1}^\pm(w)$ for all $i$.
For the  relations (12.4) we are thus reduced to prove that 
for any $\jb=(j_1,j_2,...,j_{m+1})\in\Pc_{m+1}^{no,k}$,
if $\ib=(j_1,j_2,...,j_m)$ and $\oj_{m+1}=n$ 
then
$$\eb_{i,-1} (\wedge u_\jb)=
\eb_{i,-1}(\wedge u_\ib)\wedge u_{j_{m+1}},\qquad 
\fb_{i,1} (\wedge u_\jb)=
\fb_{i,1}(\wedge u_\ib)\wedge u_{j_{m+1}}, 
\leqno(12.6)$$
for all $i=1,2,...,n-2$.
Let us prove the first equality, the second beeing quite similar.
Let $P\in\CC[X_1^{\pm 1},...,X_m^{\pm 1}]$ and
$\kb\in\Pc_{m}^0$ be such that
$\wedge u_{_\ib}=\wedge P v_{_\kb}$.
In order to simplify the notations, we omit the symbol $\otimes$ when no
confusion is possible. Put $a=-\uj_{m+1}$. Then,

$$\matrix
\eb_{_{i,-1}} (\wedge u_\jb)
&=\wedge PX_{_{m+1}}^a\eb_{_{i,-1}}(v_{_\kb} v_{_{n}}),\hfill\cr\cr
&=\wedge PX_{_{m+1}}^a\eb_{_{i,-1}}(v_{_\kb})v_{_{n}},\hfill
\endmatrix
$$

\noindent from the formula in proposition 10. Then, use
(12.3) to conclude. As for the relations (12.5), we have to prove
the following formula :
for any $\jb=(j_1,j_2,...,j_{m+2})\in\Pc_{m+2}^{no,k}$ 
set $\ib=(j_1,j_2,...,j_m)$. 
If $\oj_{m+1}=n-1$, $\oj_{m+2}=n$, and $\uj_{m+1}=\uj_{m+2}$, then

$$
\eb_{n-1,-1} (\wedge u_\jb)=
\eb_{n-1,-1}(\wedge u_\ib)\wedge u_{j_{m+1}}\wedge u_{j_{m+2}}, 
\leqno(12.7)$$
$$\fb_{n-1,1} (\wedge u_\jb)=
\fb_{n-1,1}(\wedge u_\ib)\wedge u_{j_{m+1}}\wedge u_{j_{m+2}}. 
\leqno(12.8)$$
Let us prove (12.8) for instance. 
Let $P\in\CC[X_1^{\pm 1},...,X_m^{\pm 1}]$, 
$k,l\in\NN$ and $\kb\in\Pc_{m-l-k}^0$ be such that
$$\wedge u_{_\ib}=\wedge P v_{_\kb} v_{_{n-1}}^l v_{_n}^k\quad{\and}\quad
\kb^{-1}\{n-1,n\}=\emptyset.$$
Put $a=-\uj_{m+1}=-\uj_{m+2}$. Then,

$$\matrix
&\fb_{_{n-1, 1}} (\wedge u_\jb)&
=\wedge PX_{_{m+1}}^aX_{_{m+2}}^a\fb_{_{n-1,1}}
(v_{_\kb} v_{_{n-1}}^l v_{_n}^k v_{_{n-1}} v_{_n}),
\hfill\cr\cr
(12.9)\hskip1cm\hfill&&=\wedge PX_{_{m+1}}^aX_{_{m+2}}^aT_{_{m,m-k+1}}
\fb_{_{n-1, 1}}(v_{_\kb} v_{_{n-1}}^{l+1} v_{_n}^{k+1}),
\hfill\cr\cr
(12.10)\hfill&&=\wedge PX_{_{m+1}}^aX_{_{m+2}}^aT_{_{m,m-k+1}}
v_{_\kb}\fb_{_{n-1,1}}(v_{_{n-1}}^{l+1}) v_{_n}^{k+1},
\hfill\cr\cr
(12.11)\hfill&&=\wedge PX_{_{m+1}}^aX_{_{m+2}}^aT_{_{m,m-k+1}}
v_{_\kb}\fb_{_{n-1,1}}(v_{_{n-1}}^l) v_{_{n-1}} v_{_n}^{k+1},
\hfill\cr\cr
(12.12)\hfill&&=\wedge PX_{_{m+1}}^aX_{_{m+2}}^a
\fb_{_{n-1,1}}(v_{_\kb} v_{_{n-1}}^l v_{_n}^k) v_{_{n-1}} v_{_n},
\hfill\cr\cr
(12.13)\hfill&&=\fb_{_{n-1,1}}(\wedge u_\ib)\wedge u_{j_{m+1}}\wedge 
u_{j_{m+2}}.\hfill
\endmatrix$$ 
The equalities (12.9), (12.10) and (12.12) are immediate from
the formulas in propositions 7 and 10,
while  the equality (12.13) follows from 
(12.3). As for the equality (12.11),  first note that  

$$\matrix
\fb_{_{n-1,1}}(v_{_{n-1}}^{l+1})&=
q^{-n}p^{1-n/n}
\sum_{k=1}^{l+1}q^{1-k}T_{_{k,l}}Y_{_{l+1}}^{-1} v_{_{n-1}}^l v_{_n},
\hfill\cr\cr
&=q^{-l-n}p^{1-n/n}Y_{_{l+1}}^{-1}v_{_{n-1}}^lv_{_n}+
q^{-n}p^{1-n/n}\sum_{k=1}^lq^{1-k}T_{_{k,l}}Y_{_{l+1}}^{-1}v_{_{n-1}}^lv_{_n},
\hfill\cr\cr
&=q^{-l-n}p^{1-n/n}Y_{_{l+1}}^{-1}v_{_{n-1}}^lv_{_n}+
q^{-n}p^{1-n/n}\sum_{k=1}^lq^{1-k}T_{_{k,l-1}}Y_{_l}^{-1}v_{_{n-1}}^{l-1}v_{_n}
v_{_{n-1}}+ \hfill\cr\cr
&\hskip2cm+q^{-n}p^{1-n/n}\sum_{k=1}^lq^{1-k}(q-q^{-1})T_{_{k,l-1}}
Y_{_{l+1}}^{-1}v_{_{n-1}}^lv_{_n},\hfill\cr\cr
&=Av_{_n}+\fb_{_{n-1,1}}(v_{_{n-1}}^l)v_{_{n-1}},\hfill
\endmatrix$$
where $A$ is an expression in the first $l$ components.
In this computation we have used the proposition 10, 
the equalities
$$T_l=T_l^{-1}+q-q^{-1}\quad\and\quad T_l^{-1}Y_{l+1}^{-1}=Y_l^{-1}T_l,$$
and the formulas for the $T_i$'s given in proposition 7.
Now, in order to get (12.11), it is enough to observe that 
$$\wedge PX_{_{m+1}}^aX_{_{m+2}}^aT_{_{m,m-k+1}}v_{_\kb}Av_{_n}^{k+2}$$ 
vanishes (see [KMS, Lemma 2.2.]). 
\qed

\vskip3mm

\noindent{\bf 13.} Let us now consider the classical
case, i.e. $q=1$. Let $\Ab=\CC[z^{\pm 1},D^{\pm 1}]$
be the algebra of polynomial difference operators in one variable $z$.
Suppose that $p\in\CC^\times$ is generic.
Let us recall that $z$ and $D$ satisfy the commutation relation
$D\,z=p\,z\,D$.
The algebra of matrices with coefficients in $\Ab$ is
denoted by $\gln(\Ab)$. It may be viewed as a Lie algebra
with the usual commutator. Let $\sln(\Ab)\subset\gln(\Ab)$
be the derived Lie subalgebra, i.e. $\sln(\Ab)=[\gln(\Ab),\gln(\Ab)]$.
It is known that $\sln(\Ab)\subset\gln(\Ab)$ is the subset of matrices
with trace in $[\Ab,\Ab]$.
Since $\sln(\Ab)$ is perfect, it admits a universal central extension
(see [G] and [KL] for more details). Let $\sltord$ be the 
complex Lie algebra generated by 
$\eb_{i,k}$, $\fb_{i,k}$, $\hb_{i,k}$,
where $i=0,1,...,n-1$, $k\in\ZZ$, and a central element
$\cb$, modulo the relations  
$$[\hb_{i,k},\hb_{j,l}]=d^{km_{ij}}k\delta(k=-l) a_{ij}\cb,$$
$$[\hb_{i,k},\eb_{j,l}]=d^{km_{ij}}a_{ij}\eb_{j,k+l},\quad
[\hb_{i,k},\fb_{j,l}]=-d^{km_{ij}}a_{ij}\fb_{j,k+l},$$
$$d^{-m_{ij}}[\eb_{i,k+1},\eb_{j,l}]-[\eb_{i,k},\eb_{j,l+1}]=
d^{-m_{ij}}[\fb_{i,k+1},\fb_{j,l}]-[\fb_{i,k},\fb_{j,l+1}]=0,$$
$$[\eb_{i,k},\fb_{j,l}]=\delta(i=j)(\hb_{i,k+l}+k\delta(k=-l)\cb),$$
$$\ad_{\eb_{i,0}}^{^{1-a_{ij}}}(\eb_{j,k})=
\ad_{\fb_{i,0}}^{^{1-a_{ij}}}(\fb_{j,k})=0\quad{\i\f}\quad i\neq j.$$
The algebra $\Utor_{|q=1}$ is the enveloping algebra of $\sltord$.
It is proved in [MRY] that if $d=1$ then $\sltord$ is 
isomorphic to the universal central extension, denoted $\sltor$,
of $\sln[x^{\pm 1},y^{\pm 1}]$. It is proved in [K] that
this Lie algebra can be described 
as follows : set $\Bb=\CC[x^{\pm 1},y^{\pm 1}]$ and 
$\Omega_{\Bb}=\Bb\,\d x\oplus \Bb\,\d y$,
then $\sltor=(\sln\otimes\Bb)\oplus(\Omega_{\Bb}/\d\Bb),$
with the bracket such that $\Omega_{\Bb}/\d\Bb$ is central and
$$[a\otimes f,b\otimes g]=[a,b]\otimes fg+(a|b)\,(\d f)g,
\qquad\forall a,b\in\sln,\quad\forall f,g\in\Bb,$$
where $(\cdot|\cdot)$ is the normalized Killing form of $\sln$, i.e.
$(a|b)$ is the trace of $ab$, and $\d f$ is
the differential of $f$.
The following result explains why toroidal
algebras can be represented by difference operators,
for instance as in the previous sections.
As usual, we denote by $E_{ab}$, $a,b=1,2,...,n$,
the elementary matrices of $\gln$.

\vskip3mm

\noindent{\bf Theorem 13.1.} {\it Set $d=p^{1/n}$. The map 
$$\matrix
\cb\mapsto 0,\hfill&&\cr\cr
\eb_{i,k}\mapsto p^{ki/n}E_{i,i+1}\otimes D^{-k},\quad\hfill&
\eb_{0,0}\mapsto E_{n1}\otimes z,\hfill\cr\cr
\fb_{i,k}\mapsto p^{ki/n}E_{i+1,i}\otimes D^{-k},\quad\hfill&
\fb_{0,0}\mapsto E_{1n}\otimes z^{-1},\hfill\cr\cr
\hb_{i,k}\mapsto p^{ki/n}(E_{ii}-E_{i+1,i+1})\otimes D^{-k},\quad\hfill&
\hb_{0,0}\mapsto (E_{nn}-E_{11})\otimes 1,\hfill
\endmatrix$$
where $k\in\ZZ$ and $i\neq 0$,
extends uniquely to a Lie algebra homomorphism $\pi\,:\,\sltord\to\sln(\Ab)$
such that $(\sltord,\pi)$ is the universal central extension of
$\sln(\Ab)$.
}

\vskip3mm

\noindent{\it Proof.} Let $g(z)$ be the invertible element of $\gln(\Ab)$
defined as
$$g(z)(v_i)=z^{-\delta(i=n)}v_{\overline{1+i}},\quad\forall i=1,2,...,n.$$
The conjugation by $g(z)$ is an automorphism of the associative algebra
$\gln(\Ab)$ preserving the Lie subalgebra $\sln(\Ab)$. As for the
$\eb_{i,k}$'s a direct computation gives
$$\matrix
g(z)\circ(E_{i,i+1}\otimes D^{-k})\circ g(z)^{-1}=E_{i+1,i+2}\otimes D^{-k},
\cr\cr
g(z)^2\circ(E_{n-1,n}\otimes D^{-k})\circ g(z)^{-2}=
p^{-k}E_{12}\otimes D^{-k},
\endmatrix$$
if $i=1,2,...,n-2$. Similar formulas hold for the images of 
$\fb_{i,k}$ and $\hb_{i,k}$. The argument in the beginning of the proof 
of the theorem 12 implies that the map $\pi$ such that
$$\matrix
\eb_{i,k}\mapsto p^{ki/n}E_{i,i+1}\otimes D^{-k},\hfill\cr\cr
\fb_{i,k}\mapsto p^{ki/n}E_{i+1,i}\otimes D^{-k},\hfill\cr\cr
\hb_{i,k}\mapsto p^{ki/n}(E_{ii}-E_{i+1,i+1})\otimes D^{-k},\hfill
\endmatrix$$
if $k\in\ZZ,\,i\neq 0$, and
$$\matrix
\eb_{0,k}\mapsto p^kg(z)(E_{n-1,n}\otimes D^{-k})g(z)^{-1}=
E_{n1}\otimes zD^{-k},\hfill\cr\cr
\fb_{0,k}\mapsto p^kg(z)(E_{n,n-1}\otimes D^{-k})g(z)^{-1}=
E_{1n}\otimes D^{-k}z^{-1},\hfill\cr\cr
\hb_{0,k}\mapsto p^kg(z)((E_{n-1,n-1}-E_{nn})\otimes D^{-k})g(z)^{-1}=
(p^kE_{nn}-E_{11})\otimes D^{-k},\hfill
\endmatrix$$
is a morphism of Lie algebras $\sltord\to\sln(\Ab)$. 
As for the surjectivity of $\pi$ let first note that
$$[\Ab,\Ab]=\bigoplus_{(l,k)\neq (0,0)}\CC\cdot z^kD^l.$$ 
Then, as a vector space, $\sln(\Ab)$ is the sum of $\sln\otimes\Ab$
and the subspace of $\gln(\Ab)$ of diagonal matrices with coefficients in
$[\Ab,\Ab]$. The horizontal Lie algebra (i.e. the Lie subalgebra
generated by $\eb_{i,0},\fb_{i,0},\hb_{i,0}$ with $i=0,1,...,n-1$)
and the vertical Lie algebra (generated by $\eb_{i,k},\fb_{i,k},\hb_{i,k}$
with $i\not= 0$) are isomorphic to $\slh$, the affine Lie algebra of 
type $A^{(1)}_{n-1}$ :  
the projection of $\slh$ onto each of them preserves the $\ZZ$-gradation
and thus the kernel is trivial. As a consequence,
the elements $E_{ij}\otimes z^k$, $E_{ij}\otimes D^l,$
with $i\not= j$,
and $(E_{ii}-E_{i+1i+1})\otimes z^k$, $(E_{ii}-E_{i+1i+1})\otimes D^l$
are in the image of $\pi$. Moreover, since  $p$ is generic, we have 
$$z^kD^l={1\over 1-p^{kl}}[z^k,D^l]\qquad k,l\not= 0$$
and so $E_{ij}\otimes z^kD^l,$ where $i\not= j$,
and $(E_{ii}-E_{i+1i+1})\otimes z^kD^l$
belong to the image of $\pi$ too. Thus we need only to prove that 
$E_{ii}\otimes z^kD^l\in\Im(\pi)$, for $(k,l)\not= (0,0)$.
In the case $k,l\not=0$ this follows  from 
$$[E_{ii+1}\otimes z^k,E_{i+1i}\otimes D^l]=z^kD^l(E_{ii}-p^{kl}E_{i+1i+1}),$$
$$[E_{ii+1}\otimes D^l,E_{i+1i}\otimes z^k]=z^kD^l(p^{kl}E_{ii}-E_{i+1i+1}),$$
and from the fact that $p$ is not a root of unity. In the other cases
use :
$$z^k={1\over 1-p^{-k}}[D,D^{-1}z^k],\qquad 
D^l={1\over 1-p^l}[z,z^{-1}D^l]. $$
Both algebras $\sltord$ and $\gln(\Ab)$ are graded by 
$\ZZ\times Q$, where $Q$ is the 
root lattice of $\slh$.  
Set 
$$\matrix
\deg(\eb_{i,k})=(k,\a_i),\qquad\hfill&\deg(\fb_{i,k})=(k,-\a_i),\hfill\cr\cr
\deg(\hb_{i,k})=(k,0),\qquad\hfill&\deg(\cb)=(0,0),\hfill
\endmatrix$$
and
$$\deg(E_{j,j+1}\otimes z^lD^{-k})=(k,l\delta+\a_j),\qquad
\deg(E_{j+1,j}\otimes z^lD^{-k})=(k,l\delta-\a_j),$$
$$\deg(E_{jj}\otimes z^lD^{-k})=(k,l\delta),$$
where $\a_0,\a_1,...,\a_{n-1}$ are the simple roots of $\slh$,
$\delta=\a_0+\a_1+\cdots+\a_{n-1}$ and $i=0,1,...,n-1$.
The map $\pi$ is graded. The subspace of $\sln(\Ab)$ of degree
$(k,\alpha)$ is one-dimensional if $\a$ is a real root of $\slh$ and
zero-dimensional if $\a$ is non-zero and is not a root of $\slh$.
On the other hand the Lie algebra $\sltord$ may be viewed as an 
integrable module over the horizontal Lie subalgebra 
which is isomorphic to $\slh$. Thus one can prove
as in [MRY] that the subspace of $\sltord$ of degree $(k,\alpha)$ 
is one-dimensional if $\a$ is a real root of $\slh$ and
zero-dimensional if $\a$ is non-zero and is not a root of $\slh$.
It follows that the degree
of an element of $\Ker\pi$ is in $\ZZ\times(\ZZ\cdot\delta)$
and that $(\sltord,\pi)$ is a central extension of $\sln(\Ab)$.
The Lie algebra $\sltord$ is perfect :
it is a direct consequence of the defining relations of $\sltord$.
Thus, to prove that $(\sltord,\pi)$ is universal it is enough to show that any 
central extension $(\ten,\rho)$ of $\sltord$ splits. 
Consider such an extension. To any node
$i=0,1,...,n-1$ of the Dynkin diagram of $\slh$ we associate a vertical
Lie subalgebra, $\sen_i$,  of $\sltord$ by removing the generators $\eb_{i,k}$,
$\fb_{i,k}$ and $\hb_{i,k}$. Since the Lie algebras  $\sen_i$
are isomorphic to $\slh$,  the restriction of $(\ten,\rho)$ to each of the
$\sen_i$'s splits. For each $i$ choose such a splitting
$f_i:\sen_i\rightarrow \rho^{-1}(\sen_i)$.
If $i$ and $j$ are distincts, the Lie algebra $\sen_i\cap\sen_j$ is the direct
sum of two affine Lie algebras. Thus there exists at most one morphism from 
$\sen_i\cap\sen_j$ to $\rho^{-1}(\sen_i\cap\sen_j)$
(see [G, Lemma 1.5]). As a consequence  the splittings
$f_i$'s  glue together in a splitting 
of $(\ten,\rho)$ and  we are done. 
\qed

\vskip3mm

\noindent A similar theorem holds in the case $p\to 1$. Let
$\sltorl$ be the complex Lie algebra generated by  
$\eb_{i,k}$, $\fb_{i,k}$, $\hb_{i,k}$,
where $i=0,1,...,n-1$, $k\in\NN$, modulo the relations  
$$[\hb_i(z),\hb_j(w)]=0,$$
$$[\hb_i(z),\eb_j(w)]={a_{ij}\over w-z+m_{ij}}(\eb_j(z)-\eb_j(w)),$$
$$[\hb_i(z),\fb_j(w)]=-{a_{ij}\over w-z+m_{ij}}(\fb_j(z)-\fb_j(w)),$$
$$(z-w-m_{ij})[\eb_i(z),\eb_j(w)]=
(z-w-m_{ij})[\fb_i(z),\fb_j(w)]=0,$$
$$[\eb_i(z),\fb_j(w)]={\delta(i=j)\over w-z}(\hb_i(z)-\hb_i(w)),$$
$$\ad_{\eb_{i,0}}^{^{1-a_{ij}}}(\eb_j(z))=
\ad_{\fb_{i,0}}^{^{1-a_{ij}}}(\fb_j(z))=0\quad{\i\f}\quad i\neq j,$$
where $\eb_i(z)=\sum_{k\in\NN}\eb_{i,k}z^{-k-1}$,
$\fb_i(z)=\sum_{k\in\NN}\fb_{i,k}z^{-k-1}$
and $\hb_i(z)=\sum_{k\in\NN}\hb_{i,k}z^{-k-1}$.
Put $\partial=z{d\over dz}$.

\vskip3mm

\noindent{\bf Theorem 13.2.} {\it The map 
$$\matrix
\eb_{i,k}\mapsto E_{i,i+1}\otimes (\partial-i/n)^k,\quad\hfill&
\eb_{0,0}\mapsto E_{n1}\otimes z,\hfill\cr\cr
\fb_{i,k}\mapsto E_{i+1,i}\otimes (\partial-i/n)^k,\quad\hfill&
\fb_{0,0}\mapsto E_{1n}\otimes z^{-1},\hfill\cr\cr
\hb_{i,k}\mapsto (E_{ii}-E_{i+1,i+1})\otimes (\partial-i/n)^k,\quad\hfill&
\hb_{0,0}\mapsto (E_{nn}-E_{11})\otimes 1,\hfill
\endmatrix$$
where $k\in\ZZ$ and $i\neq 0$,
extends uniquely to a Lie algebra homomorphism 
$\pi\,:\,\sltorl\to\sln(\CC[z^{\pm 1},\partial])$
such that $(\sltorl,\pi)$ is the universal central extension of
$\sln(\CC[z^{\pm 1},\partial])$.
}

\vskip3mm

\noindent{\bf Remark 13.} The algebras $\sltord$ and $\sltorl$
admit a presentation similar to the double-loop presentation
of $\sltor$. Let us recall it. Fix a complex unital associative
algebra $\Ab$. Set $Id=\sum_{i=1}^nE_{ii}\in\gln$ and
$$\matrix
[a,b]_+=ab+ba-{2\over n}(a|b)\,Id,\qquad\hfill&\forall a,b\in\sln,\cr\cr
[f,g]_+=fg+gf,\qquad\hfill&\forall f,g\in\Ab.
\endmatrix$$
Let $\Ib\subset\Ab\otimes\Ab$ be the linear span of the elements
$$f\otimes g+g\otimes f\quad{\and}\quad 
fg\otimes h-f\otimes gh-g\otimes hf$$
for all $f,g,h\in\Ab$. Denote by 
$\langle\cdot|\cdot\rangle\,:\,
\Ab\otimes\Ab\rightarrow\Ab\otimes\Ab/\Ib$
the projection. The first cyclic homology group $HC_1(\Ab)$ is
the kernel of the map 
$$\langle\Ab|\Ab\rangle\longrightarrow [\Ab,\Ab],\quad
\langle f|g\rangle\mapsto [f,g].$$
As a vector space, $\sln(\Ab)$ is the direct sum of
$\sln\otimes\Ab$ and $Id\otimes[\Ab,\Ab]$.
The bracket on $\sln(\Ab)$ is such that
$$\matrix
[a\otimes f, b\otimes g]={1\over n}(a|b)Id\otimes [f,g]+
{1\over 2}\,[a,b]\otimes [f,g]_++
{1\over 2}\,[a,b]_+\otimes [f,g],\hfill\cr\cr
[Id\otimes f,a\otimes g]=[a\otimes f,Id\otimes g]=a\otimes [f,g],\hfill
\endmatrix$$
where $a,b\in\sln$ and $f,g\in\Ab$. Similarly the universal central extension
of $\sln(\Ab)$ is the direct sum of $\sln\otimes\Ab$ 
and $\langle\Ab|\Ab\rangle$ with the bracket 
$$\matrix
[a\otimes f, b\otimes g]=
{1\over n}(a|b)\langle f|g\rangle+
{1\over 2}\,[a,b]\otimes [f,g]_++
{1\over 2}\,[a,b]_+\otimes [f,g],\hfill\cr\cr
[\langle f|g\rangle, \langle f'|g'\rangle]=
\langle\,[f,g]\,|\,[f',g']\,\rangle,\hfill\cr\cr
[\langle f|g\rangle, a\otimes h]=
a\otimes[[f,g]\,,\,h].\hfill\cr\cr
\endmatrix$$
In particular, the center is isomorphic to $HC_1(\Ab)$ (see [BGK, 359-360]
for more details).

\vskip3mm

\noindent{\bf Appendix.}
Recall that $q$ is a prime power
and $\FF$ is the field with $q^2$ elements.
Denote by $\KK=\FF((z))$ the field of Laurent power series
and by $\Bc^n$ the set of $n$-steps periodic flags (see section 6).
Let $\CC_{GL_m(\KK)}[\Bc^n\times\Bc^n]$ 
be the convolution algebra of invariant complex functions 
supported on a finite number of $GL_m(\KK)$-orbits,
where the convolution product
$$\star\,:\,\CC_{GL_m(\KK)}[\Bc^n\times\Bc^n]\otimes\CC_{GL_m(\KK)}
[\Bc^n\times\Bc^n]
\rightarrow\CC_{GL_m(\KK)}[\Bc^n\times\Bc^n]$$
is defined as
$$f\star g\, (L'',L)=\sum_{L'\in\Bc^n}f(L'',L')\cdot g(L',L).$$
If $i=0,1,...,n-1$ let 
$m_i,\chi^\pm_i,\chi^0\in\CC_{GL_m(\KK)}[\Bc^n\times\Bc^n]$ be such that 

\vskip2mm

\itemitem{\bf .} $m_i(L',L)=\dim(L_i/L_0),\qquad\forall L,L'\in\Bc^n$,

\vskip2mm

\itemitem{\bf .} $\chi^\pm_i$ is the characteristic function of the set
$$\{(L^\pm,L^\mp)\in\Bc^n\times\Bc^n\,|\,L_j^-\subset L_j^+\quad\and
\quad \dim(L_j^+/L_j^-)=\delta(\oi=\oj),\quad\forall j\in\ZZ\},$$

\vskip2mm

\itemitem{\bf .} $\chi^0$ is the characteristic function of the diagonal 
in $\Bc^n\times\Bc^n$.

\vskip2mm

\noindent{\bf Proposition.} {\it
The map
$$\eb_i\mapsto q^{m_{i-1}-m_i}\chi^+_i,\quad
\fb_i\mapsto q^{m_i-m_{i+1}}\chi^-_i,\quad
\kb_i\mapsto q^{2m_i-m_{i-1}-m_{i+1}}\chi^0,$$
extends to an algebra homomorphism 
$\Uaf\to\CC_{GL_m(\KK)}[\Bc^n\times\Bc^n]$.}

\vskip3mm

\noindent{\it Proof.}
We have to prove the $q$-deformed Kac-Moody relations written in section 3.
The relations

$$\matrix
\kb_i\cdot\kb_i^{\pm 1}=1,\qquad\hfill
&\kb_i\cdot\kb_j=\kb_j\cdot\kb_i,\hfill\cr\cr
\kb_i\cdot\eb_j=q^{a_{ij}}\,\eb_j\cdot\kb_i,\qquad\hfill
&\kb_i\cdot\fb_j=q^{-a_{ij}}\,\fb_j\cdot\kb_i,\hfill
\endmatrix$$

\noindent are immediate. As for 

$$[\eb_i,\fb_j]=\delta(i=j){\kb_i-\kb_i^{-1}\over q-q^{-1}},$$

\noindent let first remark that

$$\matrix
q^{m_{i-1}-m_i}\chi^+_i\star q^{m_j-m_{j+1}}\chi^-_j=
q^{m_{i-1}-m_i+m_j-m_{j+1}+\delta(i=j)-\delta(i=j+1)}\chi^+_i\star\chi^-_j,
\hfill\cr\cr
q^{m_j-m_{j+1}}\chi^-_j\star q^{m_{i-1}-m_i}\chi^+_i=
q^{m_{i-1}-m_i+m_j-m_{j+1}+\delta(i=j)-\delta(i=j+1)}\chi^-_j\star\chi^+_i.
\hfill
\endmatrix$$

\noindent If $i\neq j$ then $\chi^+_i\star\chi^-_j=\chi^-_j\star\chi^+_i$
is the characteristic function of the set of pairs $(L',L)$ such that
$$L'_k=L_k\quad \i\f\quad \ok\neq \oj,\oi,
\quad L'_j\subset L_j,\quad L_i\subset L'_i,
\quad\dim(L_j/L'_j)=\dim(L'_i/L_i)=1.$$
Thus,
$$\matrix
[\eb_i,\fb_j]
&=\delta(i=j)q^{m_{i-1}-m_{i+1}+1}(\chi^+_i\star\chi^-_i-\chi^-_i\star\chi^+_i),
\hfill\cr\cr
&=\delta(i=j)q^{m_{i-1}-m_{i+1}+1}
(q^{2(m_i-m_{i-1})}-q^{2(m_{i+1}-m_i)})(q^2-1)^{-1}\chi^0,
\hfill
\endmatrix$$
where the last equality simply comes from 
$\sharp(\FF\PP^k)=1+q^2+...+q^{2k}$.
Since the $\eb_i$, $\fb_i$ are locally nilpotent and since the $\kb_i$
are semisimple, the Serre relations follow from general theory of $\Uaf$.
\qed

\vskip5mm

\noindent{\it Acknowledgements.}
{\eightpoint{The authors are grateful to V. Ginzburg for stimulating
discussions.}}

\vskip5mm

\centerline{\bf References}

\vskip3mm

\hangindent3cm[B]\quad\qquad\  Beck, J.: 
Braid group action and quantum affine algebras.
{\sl Comm. Math. Phys.}, {\bf 165} (1994), 555-568.

\vskip1mm

\hangindent3cm[BGK]\quad\quad Berman, S., Gao, Y., Krylyuk, Y.S.:
Quantum tori and the structure of elliptic quasi-simple Lie algebras.
{\sl J. Funct. Anal.}, {\bf 135} (1996), 339-389.

\vskip1mm

\hangindent3cm[BLM]\quad\quad Beilinson, A., Lusztig, G., MacPherson, R.:
A geometric setting for quantum groups.
{\sl Duke Math. J.}, {\bf 61} (1990), 655-675.

\vskip1mm

\hangindent3cm[C1]\quad\qquad\  Cherednik, I.: 
Double affine Hecke algebras, Knizhnik-Zamolodchikov equations, and
Macdonald's operators. 
{\sl Int. Math. Res. Notices}, {\bf 6} (1992), 171-179.

\vskip1mm

\hangindent3cm[C2]\quad\qquad\  Cherednik, I.: Induced representations 
of double affine Hecke algebras and applications.  {\sl Math. Res.
Lett.}, {\bf 1} (1994), 319-337.

\vskip1mm

\hangindent3cm[CP]\quad\qquad Chari, V., Pressley, A.:
Quantum affine algebras and affine Hecke algebras.
{\sl qalg-preprint}, {\bf 9501003}.

\vskip1mm
\hangindent3cm[D]\quad\qquad\  Drinfeld, V.: 
A new realization of Yangians and quantized affine algebras.
{\sl Soviet. Math. Dokl.}, {\bf 36} (1988), 212-216.

\vskip1mm

\hangindent3cm[G]\quad\qquad\  Garland, H.: The arithmetic theory of
loop groups. {\sl Publ. I.H.E.S.}, {\bf 52} (1980), 5-136.

\vskip1mm

\hangindent3cm[GG]\quad\qquad Grojnowski, I., Garland, H.:
'Affine' Hecke algebras associated to Kac-Moody groups.
{\sl Preprint}, {\bf } (1995).

\vskip1mm

\hangindent3cm[GKV]\quad\quad Ginzburg, V.,  Kapranov, M.,  Vasserot, E.:
Langlands reciprocity for algebraic surfaces. {\sl Math. Res. Lett.}, {\bf 2}
(1995), 147-160.

\vskip1mm

\hangindent3cm[GL]\quad\qquad Grojnowski, I., Lusztig, G.: On bases of 
irreducible
representations of quantum $GL_n$. {\sl Contemp. Math.}, {\bf 139} (1992).

\vskip1mm

\hangindent3cm[GRV]\quad\quad Ginzburg, V.,  Reshetikhin, N.,  Vasserot, E.:
Quantum groups and flag varieties.  {\sl Contemp. Math.}, {\bf 175} (1994),
101-130.

\vskip1mm

\hangindent3cm[GV]\quad\qquad Ginzburg, V., Vasserot, E.:
Langlands reciprocity for affine quantum groups of type $A_n$.
{\sl Internat. Math. Res. Notices}, {\bf 3} (1993), 67-85.

\vskip1mm

\hangindent3cm[H]\quad\qquad\  Hayashi, T.: $Q$-analogues of Clifford and 
Weyl algebras - spinor and oscillator representations of quantum
enveloping algebras. {\sl Comm. Math. Phys.}, {\bf 127} (1990), 129-144.

\vskip1mm

\hangindent3cm[IM]\quad\qquad\ Iwahori, N., Matsumoto, H.:
On some Bruhat decompositions and the structure of Hecke rings of $p$-adic 
Chevalley groups. {\sl Pub. I.H.E.S.}, {\bf 25} (1965), 5-48.

\vskip1mm

\hangindent3cm[JKKMP]\quad Jimbo, M., Kedem, R., Konno, H., Miwa, T., 
Petersen, J.: Level-0 structure of level-1 $U_q({{\widehat{\frak{sl}}_2}})$-
modules and Mac-Donald polynomials. {\sl J. Phys.}, {\bf A28} (1995), 5589.

\vskip1mm

\hangindent3cm[K]\quad\qquad\  Kassel, C.: 
K\"ahler differentials and coverings of complex simple Lie algebras extended
over a commutative algebra.
{\sl J. Pure Appl. Algebra}, {\bf 34} (1985), 265-275.

\vskip1mm

\hangindent3cm[KL]\quad\qquad Kassel, C., Loday, J.-L.: Extensions centrales
d'alg\`ebres de Lie. {\sl Ann. Inst. Fourier, Grenoble}, {\bf 32(4)}
(1982), 119-142.

\vskip1mm

\hangindent3cm[KMS]\quad\quad Kashiwara, M., Miwa, T., Stern, E.: Decomposition
of $q$-deformed Fock space. {\sl Selecta Mathematica, New Series}, {\bf 1}
(1995), 787.

\vskip1mm

\hangindent3cm[MRY]\quad\quad Moody, R.V., Rao, S.E., Yokonuma, T.: 
Toroidal Lie algebras and vertex representations.
{\sl Geom. Dedicata}, {\bf 35} (1990), 283-307.

\vskip1mm

\hangindent3cm[PS]\quad\qquad Pressley, A., Segal, G.: Loop groups. 
{\sl Oxford Mathematical Monographs}, 1986.

\vskip1mm

\hangindent3cm[STU]\quad\quad Saito, Y., Takemura, K., Uglov, D.:
Toroidal actions on level 1 modules of $U_q(\slh)$.
{\sl qalg-preprint}, {\bf 9702024}.

\hangindent3cm[TU]\quad\qquad Takemura, K.,  Uglov, D.: Level-0 action
of $U_q(\widehat{\frak sl}_n)$ on the $q$-deformed Fock spaces.
{\sl qalg-preprint}, {\bf 9607031}.

\vskip1mm

\hangindent3cm[VV]\quad\qquad Varagnolo, M.,  Vasserot, E.: 
Schur duality in the toroidal setting. {\sl Comm. Math. Phys.}, {\bf 182}
(1996), 469-484.

\vskip3cm
{\eightpoint{
$$\matrix\format\l&\l&\l&\l\\
\phantom{.} & {\text{Michela Varagnolo}}\phantom{xxxxxxxxxxxxx} &
{\text{Eric Vasserot}}\\
\phantom{.}&{\text{Dipartimento di Matematica}}\phantom{xxxxxxxxxxxxx} &
{\text{D\'epartement de Math\'ematiques}}\\
\phantom{.}&{\text{Universit\`a di Tor Vergata}}\phantom{xxxxxxxxxxxxx} &
{\text{Universit\'e de Cergy-Pontoise}}\\
\phantom{.}&{\text{via della Ricerca Scientifica}}\phantom{xxxxxxxxxxxxx} & 
{\text{2 Av. A. Chauvin}}\\
\phantom{.}&{\text{00133 Roma}}\phantom{xxxxxxxxxxxxx} & 
{\text{95302 Cergy-Pontoise Cedex}}\\
\phantom{.}&{\text{Italy}}\phantom{xxxxxxxxxxxxx} & 
{\roman{France}}\\
&{\text{email: varagnol\@axp.mat.utovrm.it}}\phantom{xxxxxxxxxxxxx} &
{\text{email: vasserot\@math.pst.u-cergy.fr}}
\endmatrix$$
}}

\enddocument